\documentclass[12pt]{article}

\usepackage{a4wide}
\usepackage{amssymb}
\usepackage{amsfonts}
\usepackage{amsmath}
\input xy
\xyoption{arrow} \xyoption{matrix}

\date{}

\newtheorem{proposition}{Proposition}[section]
\newtheorem{theorem}[proposition]{Theorem}
\newtheorem{lemma}[proposition]{Lemma}

\newtheorem{corollary}[proposition]{Corollary}

\def\Hom{{\rm Hom}}
\def\der{\partial }

\def\nFM0{{\nu }_{F,M_0}}
\def\nFN0{{\nu }_{F,N_0}}
\def\nGN0{{\nu }_{G,N_0}}

\def\N0{ {\bf N}_0 }

\def\t{\otimes}

\def\ra{\rightarrow}

\def\Xpm{X^{\pm }}

\def\s{\sigma}

\def\l1{{\lambda}_1}

\def\a{\alpha}
\def\a0{ {\alpha }_0}
\def\a1{ {\alpha }_1}

\def\l{\lambda}

\def\nFGM0{{\nu }_{F,G,M_0}}

\def\nFN0{{\nu}_{F,N_0}}

\def\sm{{\sigma}^m}

\def\sm1{{\sigma}^{-1}}

\def\smtp1{{\sigma}^{-t+1}}

\def\S1{S^{-1}}

\def\Xpm1{X^{\pm 1}_1}

\def\sPM1{{\sigma }^{\pm 1}}
\def\sMP1{{\sigma }^{\mp 1 }}

\def\d{\delta}

\def\di{{\rm d.ind}}

\def\L{\Lambda}

\def\Ytm1{Y^{t-1}}
\def\Yim1{Y^{i-1}}

\def\CN{{\cal N}}

\def\CF{{\cal F}}

\def\Aut{{\rm Aut}}

\def\Der{{\rm Der }}
\def\ad{{\rm ad }}
\def\dim{{\rm dim }}

\def\ker{ {\rm ker } }

\def\D{ \Delta }
\def\xx{ {\bf x} }
\def\yy{ {\bf y} }

\def\SL2Z{ {\rm SL}_2({\bf Z}) }

\def\th{ \theta }

\def\Gp1{ G^{1 , 1 } }
\def\P11{ P^{-1 , 1 } }
\def\Pp1{ P^{1 , 1 } }

\def\th{\theta}

\def\nCLsr{{}^\nu\kern-2pt {\cal L}^{\sigma , \rho  }}
\def\nP{{}^\nu \kern-2pt P}
\def\nL{{}^\nu\kern-2pt L}
\def\nLL{{}^\nu\kern-2pt \Lambda}
\def\nPsr{{}^\nu\kern-2pt P^{\sigma , \rho  }}
\def\nLsr{{}^\nu\kern-2pt L^{\sigma , \rho  }}
\def\nuCL{{}^\nu\kern-2pt  {\cal L}}
\def\nCLsr{{}^\nu\kern-2pt {\cal L}^{\sigma , \rho  }}
\def\nCL1m{{}^\nu\kern-2pt {\cal L}^{-1 , 1  }}
\def\x1nu{x^\frac{1}{\nu}}
\def\xm1nu{x^{-\frac{1}{\nu}}}

\def\CN{{\cal N}}
\def\ra{\rightarrow }

\def\CC{ {\cal C}}

\def\nAM0{{\nu }_{{\cal A},M_0}}
\def\nAN0{{\nu }_{{\cal A},N_0}}

\def\Der{ {\rm Der }}

\def\det{ {\rm det }}
\def\ad{ {\rm ad }}

\def\bx{\overline{x}}
\def\by{\overline{y}}

\def\ga{\mathfrak{a}}
\def\gb{\mathfrak{b}}

\def\gn{\mathfrak{n}}
\def\gm{\mathfrak{m}}

\def\Hom{{\rm Hom}}

\def\di!{\frac{\der^i}{i!}}
\def\dik!{\frac{\der^k_i}{k!}}
\def\bd{\overline{\delta}}

\def\Fp{\mathbb{F}_p}

\def\nsder{{\rm nsder}}
\def\sder{{\rm sder}}

\def\ID{{\rm ID}}

\begin{document}

\author{V. V. \  Bavula %(simderharp.tex)
}

\title{Simple derivations of differentiably simple Noetherian commutative rings in prime characteristic}

\maketitle
\begin{abstract}
Let $R$ be  a  differentiably simple Noetherian commutative ring
of characteristic $p>0$ (then $(R, \gm )$ is local with $n:= {\rm
emdim} (R)<\infty$). A short proof is given of the Theorem of
Harper \cite{Harper61} on classification of differentiably simple
Noetherian commutative rings in prime characteristic. The main
result of the paper is that  there  exists a  nilpotent simple
derivation $\d $ of the ring $R$ such that if $\d^{p^i}\neq 0$
then $\d^{p^i}(x_i)=1$ for some $x_i\in \gm$. The derivation $\d $
is given  explicitly and  it is unique up to the action of the
group ${\rm Aut}(R)$ of {\em ring} automorphisms of $R$. Let
$\nsder (R)$ be the set of all such derivations.  Then $\nsder
(R)\simeq {\rm Aut}(R)/{\rm Aut}(R/\gm )$. The proof is based on
{\em existence} and {\em uniqueness} of an {\em iterative}
$\d$-{\em descent} (for each $\d \in \nsder (R)$), i.e. a sequence
$\{ y^{[i]}, 0\leq i<p^n\}$ in $R$ such that $y^{[0]}:=1$,
$\d(y^{[i]})=y^{[i-1]}$ and $y^{[i]}y^{[j]}={i+j\choose
i}y^{[i+j]}$ for all $0\leq i,j<p^n$. For each $\d\in \nsder (R)$,
$\Der_{k'}(R)=\oplus_{i=0}^{n-1}R\d^{p^i}$  and $k':= \ker (\d
)\simeq R/ \gm$.



$${\bf Contents}$$
\begin{enumerate}
\item Introduction. \item Differentiably simple Noetherian
commutative rings. \item Existence and uniqueness of an iterative
$\d$-descent.\item Simple derivations of differentiably simple
Noetherian commutative rings.
\end{enumerate}
\end{abstract}

%%%%%%%%%%%%%%%%%% SECTION 1 %%%%%%%%%%%%%%%%%%%%%%%%

\section{Introduction}
Throughout, ring means an associative ring with $1$ and  $p$ is a
{\em prime} number.

Let $R$ be a commutative ring. An  additive map $\d : R\ra R$ is
called a {\em derivation} of $R$ if $\d (ab) = \d (a) b +a\d (b)$
for all $a,b\in R$. If, in addition, the ring $R$ is an algebra
over a field $k$, then a derivation $\d$ is called a $k$-{\em
derivation} provided $\d (k)=0$, i.e.,  $\d$ is a $k$-linear map.
Let $\Der (R)$ be the  $R$-module of derivations of $R$. If, in
addition, the ring $R$ is an algebra over a field $k$, then
$\Der_k(R)$ denotes the $R$-module of $k$-derivations of $R$.

 An ideal $I$ of the ring $R$ is
called a {\em differential ideal} if $\d (I)\subseteq I$ for all
$\d \in \Der (R)$. The ring $R$ is called {\em differentiably
simple} if $0$ and $R$ are the only differential ideals of $R$.

\begin{theorem}\label{Harpthm}%\marginpar{Harpthm}
({\em Harper}, \cite{Harper61}) A Noetherian commutative ring $R$
of characteristic $p>0$ is differentiably simple iff it has the
form $R=k[x_1, \ldots , x_n]/(x_1^p,\ldots , x_n^p)$ where $k$ is
a field of characteristic $p$.
\end{theorem}

In his book H. Matsumura makes  the following comment on the
 Theorem of Harper (p. 206, \cite{Ma}): ``The `if' part is easy. The
proof of the `only if' part is not easy and we refer the reader to
Harper \cite{Harper61} and Yuan \cite{Yuan64}. Recently this
theorem was used by Kimura-Niitsuma \cite{Kim-Niit82} to prove the
following theorem which has been known as Kunz' conjecture.''
 Later A. Maloo \cite{Maloo93} gave a shorter proof of Theorem
 \ref{Harpthm}.

In this paper, apart from proving several  statements equivalent
 to Theorem \ref{Harpthm} (see Theorem \ref{Harpbig}) it is shown
that  essentially Theorem \ref{Harpthm} follows from Theorem 27.3
in  \cite{Ma}.

It is well known and easy to prove (Lemma \ref{Rm30Jan06}) that
each differentiably simple Noetherian ring $R$ of characteristic
$p>0$ is a local $(R, \gm )$ $k$-algebra for some subfield $k$ of
$R$ such that $R=k+\gm$. We say that a subfield $k'$ of $R$  is a
{\em coefficient field} or a {\em complement} subfield in $R$ if
$R=k'+\gm$. Clearly, each coefficient field is isomorphic to the
residue field $R/\gm$ of $R$.  Corollary \ref{cHB10Feb06}
(together with Theorem \ref{Harpbig}.(3)-(6) and Proposition
\ref{16Feb06}) gives {\em explicitly} all the coefficient fields
of $R$ (see also Theorem \ref{b10Feb06}.(2)). It is well known
that each differentiably simple Noetherian commutative ring $R$ of
characteristic $p>0$ admits a  simple  derivation (see
\cite{Yuan64} and \cite{Maxson-Ret2004}). Let $\nsder (R)$ be the
set of all {\em nilpotent simple} derivations $\d$ of the ring $R$
such that if $\d^{p^i}\neq 0$, then $\d^{p^i}(y_i)=1$ for some
element $y_i\in \gm $ where $\gm $ is a maximal ideal of $R$. The
set $\nsder (R)$ is a nonempty set (Lemma \ref{NnsdR}, see also
Theorem \ref{b10Feb06}). An action of a group $G$ on a set $X$ is
said to be fully faithful if for some/each $x\in X$ the map $G\ra
X$, $g\mapsto gx$, is a bijection.

 For a prime
number $p$, $\Fp := \mathbb{Z}/p\mathbb{Z}$ is a finite field that
contains $p$ elements. Given a derivation $\d$ of an $\Fp$-algebra
$A$, a sequence $\{ y^{[i]}, 0\leq i<p^n\}$ of elements in $A$
(where $y^{[0]}:=1$) is called an {\bf iterative} $\d$-{\bf
descent} if
$$ \d (y^{[i]})=y^{[i-1]}, \;\;\; y^{[i]}y^{[j]}={i+j\choose i}y^{[i+j]},\;\; \; 0\leq  i,j\leq
p^n-1,$$ where  $y^{[-1]}=y^{[k]}:=0$ for all $k\geq p^n$.

Let us give a list of  the main results of the paper. {\em Let $R$
be a differentiably simple Noetherian commutative ring of
characteristic $p>0$ which is not a field, $\gm $ be its maximal
ideal, $k:=R/\gm $, $n:=\dim_k(\gm / \gm^2)\geq 1$. Let $\d \in
\nsder (R)$ and $k':= \ker (\d )$. Then}
\begin{itemize}
\item (Theorem \ref{md11Feb06}) $k'$ {\em is a coefficient field
for $R$.} \item (Theorem \ref{md11Feb06}.(1)) $\d^{p^n-1}\neq 0$
{\em  and} $\d^{p^n}=0$. \item (Proposition \ref{b17Feb06})
$\Der_{k'}(R)=\bigoplus_{i=0}^{n-1} R\d^{p^i}$. \item (Theorem
\ref{md11Feb06}.(2)) {\em There exists a {\bf unique} iterative
$\d$-descent $\{ x^{[i]}, 0\leq i <p^n\}$. Then  $x^{[i]}\in \gm$,
$1\leq i <p^n$;
$$R=\bigoplus_{i=0}^{p^n-1}k'x^{[i]}= k'\langle x_0, \ldots ,
x_{n-1}\rangle\simeq k'[x_0, \ldots , x_{n-1}]/(x_0^p, \ldots ,
x_{n-1}^p),$$ where $x_j:= x^{[p^j]}$ and $x^{[i]}:=
\prod_{k=0}^t\frac{x_k^{i_k}}{i_k!}$ where $i=\sum_{k=0}^ti_kp^k$,
$0\leq i_k<p$, and} \item (Theorem \ref{f18Jan06}) {\em the
iterative $\d$-descent $\{ x^{[i]}, 0\leq i <p^n\}$ is given
explicitly: choose elements $y_0, y_1, \ldots , y_{n-1}\in \gm$
such that $\d^{p^k}(y_k)=1$ for $k=0, \ldots , n-1$; then $x_0:=
y_0$,
$$x_1:= (-1)^{p-1}\phi_0(y_1), \;\; \phi_0 (z):=
\sum_{j=0}^{p-1}(-1)^j\frac{x_0^j}{j!}\d^j(z),$$ and then
recursively, for all} $1\leq k \leq n-2$,
$$x_{k+1}:= (-1)^{p-1}\d^{p^k-1}(\prod_{l=0}^{k-1}
\frac{x_l^{p-1}}{(p-1)!}\cdot \phi_k(y_{k+1})), \;\; \phi_k (z):=
\sum_{j=0}^{p-1}(-1)^j\frac{x_k^j}{j!}\d^{p^kj}(z).$$ \item (Lemma
\ref{Ky13Feb06}, Theorem \ref{md11Feb06}) {\em The derivation $\d$
has the unique presentation via its iterative $\d$-descent,  $\d =
\sum_{i=0}^{n-1} x^{[p^i-1]}\frac{\der }{\der x_i}$, and} $\d \in
\nsder _{k'}(R):= \{ \der \in \nsder (R)\, | \, \der (k')=0\}$.
\item (Theorem \ref{md11Feb06}) $k'=\phi (R)$ {\em where} $\phi :=
\sum_{i=0}^{p^n-1} (-1)^i x^{[i]}\d^i : R=k'\oplus \gm \ra
R=k'\oplus \gm$ {\em is the projection onto} $k'$. \item
(Corollary \ref{a19Feb06}.(1)) {\em For each coefficient field}
$l$ {\em of} $R$, {\em the action
 $${\rm Aut}_l(R)\times \nsder_l (R)\ra \nsder_l (R)
$$ defined by the rule $(\s, \der )\mapsto \s \der \s^{-1}$
 is fully faithful where ${\rm
Aut}_l(R)$ is a the group of $l$-algebra automorphisms of $R$ and}
$\nsder_l(R):= \{ \der \in \nsder (R)\, | \, \der (l)=0\}$. \item
(Corollary \ref{a19Feb06}.(2)) {\em The action ${\rm Aut}(R)\times
\nsder (R)\ra \nsder (R)$ which is given by the rule  $(\s, \d
)\mapsto \s \d \s^{-1}$ has a single orbit and,  for each $\der
\in \nsder (R)$, ${\rm Fix} (\der )\simeq {\rm Aut}(k)$, and so}
$\nsder (R)\simeq {\rm Aut}(R)/{\rm Aut} (k)$.
\end{itemize}

Note that the group ${\rm Aut}_{k'}(R)$ is easily described: any
$k'$-automorphism of the algebra $R$ is uniquely determined by $n$
``polynomials'' $\s(x_i)= \sum a_{ij} x_j+\cdots $, $0\leq i\leq
n-1$, $a_{ij}\in k'$,  with $\det (a_{ij})\neq 0$ where the three
dots mean {\em any} linear combination of monomials of degree
$\geq 2$. So, the result above  gives explicitly all the elements
of $\nsder_{k'} (R)$.

In brief, almost all the results of the paper are  consequences of
two theorems on  existence  and uniqueness of an iterative
$\d$-descent (Theorem \ref{f18Jan06} and Theorem \ref{md11Feb06})
which is given  explicitly for each $\d\in \nsder (R)$. The
importance of the iterative $\d$-descent lies in the facts that
$(i)$ the iterative $\d$-descent together with $\ker (\d )$
determines uniquely and explicitly the derivation $\d \in \nsder
(R)$ (Lemma \ref{Ky13Feb06}, Theorem \ref{md11Feb06}), $(ii)$ all
the coefficient fields are precisely the kernels of derivations
from $\nsder (R)$; and the iterative $\d$-descent $\{ x^{[i]}\}$
describes explicitly the kernel of $\d$: $\ker (\d ) =
(\sum_{i=0}^{p^n-1} (-1)^ix^{[i]}\d^i)(R)$ (Theorem
\ref{md11Feb06}).

The paper is organized as follows: in Section 1, a short proof of
the Theorem of Harper is given together with some equivalent
statements (Theorem \ref{Harpbig}). Corollary \ref{cHB10Feb06}
describes explicitly all the coefficient fields. An important
technical result, Proposition \ref{16Feb06}, is proved which
states roughly that having $n$ derivations $\d_i$ with
$\d_i(x_j)=\d_{ij}$ one can produce $n$ {\em commuting}
derivations $\d_i'$ such that $\d_i'(x_i')=\d_{ij}$ and
$\d_i'^p=0$.

In Section 3, the concept of an iterative $\d$-descent is
introduced. Theorem \ref{f18Jan06} (on existence and uniqueness of
an iterative $\d$-descent) is proved.

In Section 4, the canonical form for each derivation $\d \in
\nsder (R)$ is given via the iterative $\d$-descent (Theorem
\ref{md11Feb06}, see also Lemma \ref{Ky13Feb06}), and the
coefficient fields are explicitly described (Theorem
\ref{md11Feb06}). An important bijection is established in Theorem
\ref{b10Feb06} which is used in the proof of the canonical
bijection $\nsder (R)\simeq \Aut (R)/ \Aut (R/ \gm )$ (Corollary
\ref{a19Feb06}). Finally, it is proved that $\Der_{R^\d}(R)
=\bigoplus_{i=0}^{n-1} R\d^{p^i}$ (Proposition \ref{b17Feb06}).

%%%%%%%%%%%%%%%%%% SECTION 2 %%%%%%%%%%%%%%%%%%%%%%%%

\section{Differentiably simple Noetherian commutative rings}
In this section, a short proof of Theorem \ref{Harpthm} is given
 and some equivalent statements to Theorem \ref{Harpthm} are proved  (Theorem
 \ref{Harpbig}).
 All coefficient fields of a differentiably simple Noetherian
commutative ring of prime characteristic are found explicitly (see
the remark at the end of this section).

\begin{lemma}\label{Rm30Jan06}%\marginpar{Rm30Jan06}
Let $R$ be a differentiably simple commutative ring of
characteristic $p>0$ (i.e.,  $p^k=0$ in $R$ for some $k\geq 1$).
Then $R$ is a local $\Fp$-algebra with maximal ideal $\gm$ such
that $x^p=0$ for all $x\in \gm$, and $\gm^\infty := \cap_{i\geq 1}
\gm^i=0$. If, in addition, $R$ is a Noetherian ring, then $R=
k+\gm$ for some subfield $k$ of $R$ necessarily isomorphic to the
 residue field $R/ \gm$ of $R$.
\end{lemma}

{\it Proof}. Let $\gm $ be a maximal ideal of $R$. Then the ideal
$I:=\sum_{x\in \gm }Rx^p$ is  differential (since $\d
(x^p)=px^{p-1}\d (x)=0$ for all $\d \in \Der (R)$), and
$I\subseteq \gm$; hence $I=0$ (since $R$ is a differentiably
simple ring). Therefore, $\gm = \gn (R)$ is the only maximal ideal
of $R$ where $\gn (R)$ is the  nil radical of $R$;  that is,  $(R,
\gm )$ is a local ring.

The ring of constants $C:= \bigcap_{\d \in \Der (R)}\ker \, \d$
must be a field since $R$ is a differentiably simple ring (for any
$0\neq c\in C$, $cR$ is a differential ideal of the ring $R$;
hence $cR=R$). Therefore, $R$ is an $\Fp$-algebra. Note that
$\gm^\infty $ is a differential ideal of $R$ such that $\gm^\infty
\subseteq \gm$. Hence $\gm^\infty =0$.

If, in addition, $R$ is a Noetherian ring, then $(R, \gm )$ is a
Noetherian local ring which is obviously equicharacteristic and
complete in the $\gm$-adic topology  since $\gm^i=0$ for a large
$i$.  It is well known that {\em any equicharacteristic complete
Noetherian local  commutative ring contains a coefficient field}
(\cite{Ma}, Theorem 28.3.(ii)), and so $R= k+\gm$ for some
subfield $k$ of $R$. $\Box $

So, in dealing with  a differentiably simple commutative ring $R$
of characteristic $p>0$ there is no restriction in assuming that
it is a local $\Fp$-algebra $(R, \gm )$ with $x^p=0$ for all $x\in
\gm$. That explains why in many results of the present paper these
conditions are present from the outset (aiming at possible
application to differentiably simple rings).

Let $\Fp [ h]$ be a polynomial algebra in a variable $h$. The
factor algebra $\L := \Fp [h]/(h(h-1) \cdots (h-p+1))$ is
isomorphic to the direct product of $p$ copies of the field $\Fp$.
In more detail, $\L =\oplus_{i=0}^{p-1}\Fp \th_i$ where
$1=\th_0+\th_1+\cdots +\th_{p-1}$ is a sum of primitive orthogonal
idempotents of the algebra $\L $ (i.e.,
$\th_i\th_j=\d_{ij}\th_i$ for all $i,j\in \mathbb{Z}/p\mathbb{Z}$
where $\d_{ij}$ is the {\em Kronecker delta})   and
\begin{eqnarray*}
 \th_i &:=& \frac{h(h-1) \cdots
\widehat{(h-i)}\cdots (h-p+1)}{i(i-1)\cdots 1(-1)(-2)\cdots
(i-p+1)},
\end{eqnarray*}
where the hat over a symbol means that it is missed. Using the
facts that $-i= p-i$ and $(-1)^{p-1}=1$ in $\Fp$, it follows at
once that
\begin{eqnarray*}
 \th_i &=& \frac{h(h-1) \cdots \widehat{(h-i)}\cdots
(h-p+1)}{(p-1)!} = (-1)^{p-1}\frac{h(h-1) \cdots
\widehat{(h-i)}\cdots (h-p+1)}{(p-1)!}.
\end{eqnarray*}
These presentations of $\th_i$ will be used later in calculations.

The $\Fp$-algebra automorphism $\s \in \Aut_{\Fp} (\L )$,
$h\mapsto h-1$, permutes cyclicly the idempotents $\th_i$: $\s
(\th_i)= \th_{i+1}$. It is
evident that %\marginpar{p1thi}
\begin{equation}\label{p1thi}
\prod_{i=0}^{p-1}(1-\th_i)=1-\th_0-\th_1-\cdots - \th_{p-1}=1-1=0.
\end{equation}
Note that if $\d $ is a derivation of an  $\Fp$-algebra $A$, then
so are $\d^{p^i}$, $i\geq 0$. For $a\in A$, $(\ad \, a )(x):=
ax-xa$ is an {\em inner derivation} of $A$. It is obvious that
$\Der (A) = \Der_{\Fp} (A)$ for each $\Fp$-algebra $A$ since $\d
(1)= \d (1\cdot 1)= 2\d (1)$, i.e.,  $\d (1)=0$.

\begin{lemma}\label{dpx27Jan06}%\marginpar{dpx27Jan06}
Let $R$ be a commutative $\Fp$-algebra, $\d \in \Der (R)$, $\d
(x)=1$ for some $x\in R$ such that $x^p=0$, $h:= x\d \in {\rm
End}_{\Fp }(R)$ where $x$ is identified with the $\Fp$-linear map
$r\mapsto xr$.
\begin{enumerate}
\item The $\Fp$-subalgebra of the endomorphism algebra ${\rm
End}_{\Fp }(R)$  generated by $h$ is naturally isomorphic to the
factor algebra $\L := \Fp [h]/(h(h-1) \cdots
(h-p+1))=\bigoplus_{i=0}^{p-1}\Fp\th_i$ (as above). \item For each
$i=0, \ldots , p-1$, let $\d_i:= (1-\th_i)\d$. Then
$\d_i^p=0$.\item The map $\der := \d_{p-1}= (1-\th_{p-1})\d $ is,
in fact, the derivation $\der =\d - (-1)^{p-1}
\frac{x^{p-1}}{(p-1)!}\d^p\in \Der (R)$ that satisfies the
conditions $\der (x)=1$, $\der^p=0$, and $h=x\der $.
\end{enumerate}
\end{lemma}

{\it Proof}. $1$. In the algebra ${\rm End}_{\Fp}(R)$, we have
$\d x-x\d =\d (x)=1$, hence $xh=(h-1)x$. Using this relation, we
have
$$ 0=x^p\d^p= h(h-1)\cdots (h-p+1), $$ and so there is a natural
$\Fp$-algebra epimorphism $\L \ra \Fp \langle h\rangle$. It is, in
fact, an isomorphism.  It suffices to show that the elements $1,
h, \ldots , h^{p-1}$ are $\Fp$-linearly independent in $\Fp
\langle h\rangle$. If $r:= \l_0+\l_1 h+\cdots +\l_mh^m=0$ is a
nontrivial relation  (all $\l_i\in \Fp$ and $\l_m\neq 0$, $0\leq m
\leq p-1$), then applying $(-\ad \, x)^m$ to $r$ we have $m!\l_m
x^m=0$ in ${\rm End}_{\Fp }(R)$. Evaluating this relation at $1$,
one has the relation $m!\l_m x^m=0$ in the ring $R$, and so
$0=\d^m(\l_m x^m)=m!\l_m \neq 0$, a contradiction.

Recall that the $\Fp$-algebra automorphism $\s \in \Aut_{\Fp}(\L
)$ is defined as follows: $\s (h)=h-1$. Then $\s^{-1}(h)= h+1$.
In the algebra ${\rm End}_{\Fp}(R)$, $\d x-x\d =\d (x)=1$. For
computational reasons, it is important to stress  that this
relation is equivalent to the following four relations:
\begin{eqnarray*}
 xh=\s (h) x, & x\d =h,\\
 \d h = \s^{-1}(h)\d, & \;\; \;\;\; \; \; \d x =\s^{-1}(h).
\end{eqnarray*}

$2$. Note that each $\th_i\in \Fp [h]$ and $\d
\th_i=\s^{-1}(\th_i)\d$. Then, by the equality (\ref{p1thi}),
$$ \d_i^p= \prod_{j=0}^{p-1}\s^j (1-\th_i)\cdot \d^p=
\prod_{j=0}^{p-1} (1-\th_{i+j})\cdot \d^p=\prod_{k=0}^{p-1}
(1-\th_k)\cdot \d^p=0\cdot \d^p=0.$$

$3$. Since $x^{p-1}\d^{p-1}= x^{p-2}h\d^{p-2}=
(h-p+2)x^{p-2}\d^{p-2}=\cdots = (h-p+2) \cdots (h-1)h$, we have
\begin{eqnarray*} \der &= &(1-\th_{p-1})\d = \d
-(-1)^{p-1}\frac{h(h-1)\cdots (h-p+2)}{(p-1)!}\d = \d -
(-1)^{p-1}\frac{x^{p-1}\d^{p-1}}{(p-1)!}\d \\
&=& \d -(-1)^{p-1} \frac{x^{p-1}}{(p-1)!}\d^p\in \Der (R).
\end{eqnarray*}
Then it becomes obvious that  $\der (x)=\d (x)=1$ and $x\der = x\d
-(-1)^{p-1}\frac{x^p}{(p-1)!}\d^p =x\d = h$.
 $\Box $

\begin{theorem}\label{30Jan06}%\marginpar{30Jan06}
Let $\d $ be a derivation of a commutative  $\Fp$-algebra $R$ such
that $\d^p=0$ and $\d (x)=1$ for some element $x\in R$. Then
\begin{enumerate}
\item  (Theorem 27.3, \cite{Ma}) $R=\bigoplus_{i=0}^{p-1}R^\d x^i$
where $R^\d := \ker (\d )$. \item The map $\phi :=
\sum_{i=0}^{p-1} (-1)^i \frac{x^i}{i!}\d^i: R=R^\d \oplus (x)\ra
R=R^\d \oplus (x)$ is a projection onto the subalgebra $R^\d$;
that is,  $\phi (a+bx)= a$ for all $a\in R^\d$ and $b\in R$. \item
For any $a\in R$, $a=\sum_{i=0}^{p-1}\phi
(\frac{\d^i}{i!}(a))x^i$.
\end{enumerate}
\end{theorem}

{\it Proof}. $2$. By the very definition, the map $\phi$ is a
homomorphism of $R^\d$-modules. For each $n=1, \ldots , p-1$,
$$\phi (x^n)=\sum_{i=0}^{p-1}(-1)^i{n\choose i} x^ix^{n-i}= (1-1)^nx^n=0\cdot x^n=0.$$
Therefore, $\phi (a+bx)=a$,  which follows directly from statement
$1$.

$3$. Given $a=\sum_{i=0}^{p-1} a_ix^i\in R$ where $a_i\in R^\d$,
 by statement $2$, $a_i= \phi (\frac{\d^i}{i!}(a))$.  $\Box $

\begin{lemma}\label{d16Feb06}%\marginpar{d16Feb06}
Let $A=\bigoplus_{\alpha \in \mathbb{N}^n}A_\alpha = A_0\oplus
A_+$ be an $\mathbb{N}^n$-graded ring where $A_+:=\bigoplus_{0\neq
\alpha \in \mathbb{N}^n}A_\alpha$, and let $\d$ be a derivation of
the ring $A$. Then
\begin{enumerate}
\item for each $a\in A_0$, $\d (a)=\d_0(a)+\d_+(a)$ where
$\d_0(a)\in A_0$, $\d_+(a)\in A_+$, and  $\d_0$ is a derivation of
the ring $A_0$. \item If $\d (A_+)\subseteq A_+$ and $\d (x)=1$
for some element $x=x_0+x_+$ where $x_0\in A_0$ and $x_+\in A_+$,
then $\d_0(x_0)=1$.
\end{enumerate}
\end{lemma}

{\it Proof}. $1$. Since $A= A_0\oplus A_+$, both maps $\d_0$ and
$\d_+$ are additive, by the very definition. For any $a,b\in A_0$,
$\d (ab)= \d (a) b+a\d (b)= \d_0(a)b+a\d_0(b)+c$ for some $c\in
A_+$; hence $\d_0 (ab)= \d_0(a)b+a\d_0(b)$. This means that $\d_0$
is a derivation of the ring $A_0$.

$2$. The equality $1= \d (x)= \d_0(x_0)+\d_+(x_0)+\d(x_+)$ implies
the equality  $\d_0(x_0)=1$. $\Box $

Recall that for a local commutative ring $(R, \gm )$, a subfield
$k'$ of $R$ satisfying $R= k'+\gm$ is called a coefficient field
for $R$. For a natural number $n\geq 1$, let  $$\CN_n:=\{ \alpha =
(\alpha_1, \ldots , \alpha_n)\in \mathbb{N}^n\, | \, 0\leq
\alpha_\nu <p\}.$$

\begin{proposition}\label{16Feb06}%\marginpar{16Feb06}
Let $(R, \gm )$ be a local Noetherian commutative $\Fp$-algebra
such that $x^p=0$ for all $x\in \gm$,  $k:= R/\gm $, $n:=\dim_k(
\gm / \gm^2)\geq 1$. Let  $\d_1, \ldots , \d_n$ be derivations of
the ring $R$ such that $\d_i(x_j)=\d_{ij}$ for some elements $x_1,
\ldots , x_n\in \gm$. Then there exist commuting derivations
$\d_1', \ldots , \d_n'$ of the ring $R$ and elements $x_1', \ldots
, x_n'\in \gm$ such that $\d_i'(x_j')=\d_{ij}$, $\d_1'^p=\cdots =
\d_n'^p=0$, and $(x_1', \ldots , x_n')=\gm$. Then necessarily
$k':= \bigcap_{i=1}^n \ker \, \d_i'$ is a coefficient field of $R$
and $R= \bigoplus_{\alpha \in \CN_n} k'x^\alpha$ where $x^\alpha
:= x_1^{\alpha_1}\cdots x_n^{\alpha_n}$. In particular, $R\simeq
k'[x_1, \ldots , x_n]/(x_1^p, \ldots , x_n^p)$.
\end{proposition}

{\it Remark}.  Necessarily, $\gm =(x_1, \ldots , x_n)$ since
$n=\dim_k( \gm / \gm^2)$ and $\d_i(x_j)=\d_{ij}$.

 {\it Proof}. The idea of the proof is to use repeatedly Lemma
\ref{dpx27Jan06}.(3), Theorem \ref{30Jan06}.(1), and Lemma
\ref{d16Feb06}.

Since $x_1^p=0$ and  $\d_1(x_1)=1$, by Lemma \ref{dpx27Jan06}.(3),
one can find a derivation $\d_1'$ of the ring $R$ such that
$\d_1'(x_1)=1$ and  $\d_1'^p=0$. Let $x_1':= x_1$. Then
$(x_1')=(x_1)$. By Theorem \ref{30Jan06}.(1),
$R=\bigoplus_{\alpha_1=0}^{p-1}R^{\d_1'}x_1'^{\alpha_1}$ is a
positively graded ring.

Suppose that using the derivations $\d_1, \ldots , \d_s$ and the
elements $x_1, \ldots , x_s$ we have already found commuting
derivations $\d_1', \ldots , \d_s'$ and  elements $x_1', \ldots ,
x_s'\in \gm$ such that the following conditions hold:
$\d_1'^p=\cdots = \d_s'^p=0$; $\d_i'(x_j')=\d_{ij}$ for all
$i,j=1, \ldots , s$;  $(x_1', \ldots , x_s')=(x_1, \ldots , x_s)$
(the equality of ideals); and
$$R=\bigoplus_{\alpha \in \CN_s} k_s'x^\alpha, \;\; \;
k_s':=\bigcap_{i=1}^s \ker \, \d_i'.$$ The ring $R$ is {\em
naturally} $\mathbb{N}^s$-graded, and  $R_+:=\bigoplus_{0\neq
\alpha \in \CN_s} k_s'x^\alpha = (x_1', \ldots , x_s')$. Then
$$\d_{s+1}(R_+)= \d_{s+1}((x_1', \ldots
, x_s'))=\d_{s+1}((x_1, \ldots , x_s))\subseteq (x_1, \ldots ,
x_s)=(x_1', \ldots , x_s')= R_+.$$ Write
$x_{s+1}=x_{s+1}'+x^+_{s+1}$ for some  $x_{s+1}'\in k_s'\cap \gm$
and
 $x^+_{s+1}\in R_+\subseteq \gm$. Then $(x_1', \ldots
, x_{s+1}')=(x_1, \ldots , x_{s+1})$. Since $\d_{s+1}(x_{s+1})=1$,
by Lemma \ref{d16Feb06}, $\d_{s+1}'(x_{s+1}')=1$ for some
derivation $\d_{s+1}'$ of the ring $k_s'$. One can extend the
derivation $\d_{s+1}'$ to a derivation, say $\d_{s+1}'$, of the
ring $R$ by setting $\d_{s+1}'(x_1')=\cdots = \d_{s+1}'(x_s')=0$.
Changing (if necessary) the derivation $\d_{s+1}'$ as in Lemma
\ref{dpx27Jan06}.(3), one can assume additionally that
$\d_{s+1}'^p=0$. Then the derivations $\d_1', \ldots , \d_{s+1}'$
commute,  the elements $x_1', \ldots , x_{s+1}'\in \gm$,
$\d_1'^p=\cdots = \d_{s+1}'^p=0$, $\d_i'(x_j')=\d_{ij}$ for all
$i,j=1, \ldots , s+1$, $(x_1', \ldots , x_{s+1}')=(x_1, \ldots ,
x_{s+1})$, and, by Theorem \ref{30Jan06}.(3),
$$R=\bigoplus_{\alpha \in \CN_{s+1}} k_{s+1}'x^\alpha, \;\; \;
k_{s+1}':=\bigcap_{i=1}^{s+1} \ker \, \d_i'.$$ Now, by induction
on $s$ we have
$$R=\bigoplus_{\alpha \in \CN_n} k'x^\alpha, \;\; \;
k':=\bigcap_{i=1}^n \ker \, \d_i',$$ $(x_1', \ldots , x_n')=(x_1,
\ldots , x_n)=\gm $. Hence $k'\simeq k$, and so $k'$ is a
coefficient field for $R$. $\Box $

Let $V$ and $U$ be a finite dimensional vector spaces over a field
$k$. A $k$-bilinear map $V\times U\ra k$, $ (v,u)\mapsto vu$, is
called a pairing. It is a {\em perfect pairing} if ${\rm ann}
(V):= \{ u\in U\, | \, Vu=0\}=0$ and ${\rm ann} (U):= \{ v\in V\,
| \, vU=0\} =0$. A $k$-bilinear map $V\times U\ra k$, $
(v,u)\mapsto vu$, is a perfect pairing iff one of the equivalent
conditions holds:

$(i)$ $\dim_k(V) = \dim_k(U)$ and ${\rm ann} (V)=0$;

$(ii)$ the map $V\ra U^*:= \Hom_k(U,k)$, $v\mapsto (u\mapsto vu)$,
is a bijection;

$(iii)$ the map $V\t_kU\ra {\rm End}_k(V)$, $v\t u \mapsto
(v'\mapsto v(uv'))$, is a bijection.

Next, a short proof is given of the Theorem of Harper (Theorem
\ref{Harpbig}, $( 1\Leftrightarrow 2)$), and some equivalent
statements to the Theorem of Harper are added. Note that in the
proof $(1\Rightarrow 2 )$ we use only Theorem \ref{30Jan06} and
Lemma \ref{dpx27Jan06}.

\begin{theorem}\label{Harpbig}%\marginpar{Harpbig}
Let $(R, \gm )$ be a local Noetherian commutative $\Fp$-algebra
such that $x^p=0$ for all $x\in\gm$. Let   $k:= R/\gm $,
$n:=\dim_k( \gm / \gm^2)\geq 1$, and  $T:= \Der (R)/\gm \Der (R)$.
Then the following statements are equivalent:
\begin{enumerate}
\item $R$ is a differentiably simple ring.\item $R\simeq k[x_1,
\ldots , x_n]/(x_1^p, \ldots , x_n^p)$. \item There exist
derivations $\d_1, \ldots , \d_n$ of the ring  $R$ and  elements
$x_1, \ldots , x_n\in \gm$ such that $\det (\d_i (x_j))\not\in
\gm$. \item There exist  commuting derivations $\d_1,\ldots , \d_n
$ of $R$ and  elements $x_1,\ldots , x_n\in \gm$ such that
 $\d_1^p=\cdots = \d_n^p=0$ and $\d_i(x_j)=\d_{ij}$ for all $i,j=1, \ldots , n$. \item The $k$-bilinear map $T\times
\gm /\gm^2\ra k$ given by the rule $(\bd, \bx )\mapsto \bd (\bx )
:= \d (x)+\gm $ has ${\rm ann} (T):=\{ u\in \gm / \gm^2\, | \,
Tu=0 \} =0$ and $\dim_k(T)\geq n$ where $\bd := \d +\gm \Der
(R)\in T$ and $\bx : = x+\gm^2\in \gm / \gm^2$. Equivalently, the
$k$-bilinear map is a perfect pairing. \item The $k$-linear map
$\gm / \gm^2\t_k T\ra {\rm End}_k(\gm / \gm^2)$ given by the rule
$\bx\t \bd  \mapsto (\by \mapsto \bx  \bd (\by ) )$
 is a surjection or, equivalently, is a bijection.
 \item  $R$ is a $\Der_{k'}(R)$-simple $k'$-algebra for some/any coefficient field $k'$ of $R$ (i.e.,  $R=k'\oplus \gm$).
\end{enumerate}
\end{theorem}

{\it Proof}. The implications $7({\rm some})\Rightarrow 1$ and
$7({\rm any})\Rightarrow 7({\rm some})$ are trivial.

$(4\Rightarrow 3)$: Use the same $\d_i$ and $x_j$ and the fact
that $\det (\d_i(x_j))= \det (E)=1$ where $E$ is the identity
matrix.

$(5\Rightarrow 6)$: Suppose that statement 5 holds. Let us prove
first that the pairing of statement 5 is perfect; i.e.,  ${\rm
ann} (T)=0$ and $\dim_k(T)=n$. The first condition is given. To
prove the second it suffices to show that $\dim_k(T)\leq n$ since
$\dim_k(T)\geq n$, by the assumption. Since ${\rm ann}(T)=0$, the
$k$-linear map
$$ T\ra (\gm / \gm^2)^*:= \Hom_k(\gm / \gm^2 , k), \;\; \bd
\mapsto (\by \mapsto \bd (\by )), $$ is injective. Therefore,
$\dim_k(T)\leq \dim_k((\gm / \gm^2)^*)= \dim_k(\gm/ \gm^2) = n$,
as required. So, the pairing in statement 5 is perfect, i.e.,  the
map in statement 6 is bijective. In particular,  it is surjective.

$(6\Rightarrow 5)$: Suppose that the first part of  statement 6
holds; i.e.,  the map is surjective. First, we prove that
statement 5 holds,  which then gives the fact that the map in
statement 6 is a bijection (see $(5\Rightarrow 6)$ above).  Since
the map in statement 6 is surjective, this implies that the set
${\rm ann} (T)$ is annihilated by ${\rm End}_k(\gm /\gm^2)$, and
so $ {\rm ann}(T)=0$. This gives the first condition of statement
5. On the other hand,
$$ n^2= \dim_k({\rm End}_k(\gm
/\gm^2)) \leq \dim_k(\gm/\gm^2\t_kT) = n\cdot \dim_k(T);$$ hence
$\dim_k(T)\geq n$. This gives the second condition of statement 5.
So, statement 5 holds; hence the map in statement 6 is a bijection
as was proved above in $(5\Rightarrow 6)$.  So, statements 5 and 6
are equivalent.

$(2\Rightarrow 4)$: If $R= k[x_1, \ldots , x_n]/(x_1^p ,\ldots ,
x_n^p)$ (up to isomorphism), then the partial derivatives
$\der_1:= \frac{\der}{\der x_1}, \ldots , \der_n:=
\frac{\der}{\der x_n}\in \Der_k(R)$ and the elements $x_1, \ldots
, x_n$ satisfy the conditions of statement 4.

$(2\Rightarrow 6)$: ${\rm End}_k(\gm /\gm^2)= \bigoplus_{i,j=1}^n
k\bx_i\overline{\der}_j$ since $\dim_k({\rm End}_k(\gm
/\gm^2))=n^2$ and $\bx_i\overline{\der}_j(\bx_k)=\d_{j,k}\bx_i$
for all $i,j,k$ (i.e.,  $\bx_i\overline{\der}_j$ play the role of
the matrix units).

$(2\Rightarrow 7({\rm some}))$: Since $\Der_k(R) =
\bigoplus_{i=1}^n R\der_i$, the ring $R$ is a simple
$\Der_k(R)$-module.

$(5\Rightarrow 3)$: The pairing $T\times \gm / \gm^2\ra k$ is
perfect. Choose dual bases, say $\bd_1, \ldots , \bd_n\in T$ and
$\bx_1, \ldots , \bx_n\in \gm / \gm^2$ (i.e.,  $\bd_i (\bx_j) =
\d_{i,j}$ for all $i$ and $j$). Therefore, $\det (\d_i(x_j))\equiv
1 \mod \gm$, as required.

 It remains to prove the implications
$1\Rightarrow 2$, $3\Rightarrow 2$, and $7({\rm some})\Rightarrow
7({\rm any})$.

$(1\Rightarrow 2)$: Suppose that the algebra $R$ is differentiably
simple. Then $\gm $ is not a differential ideal of $R$, i.e.,
$\d_1 (\gm )\not\subseteq \gm$ for some derivation $\d_1$ of $R$.
Then $\d_1(x_1)\not\in \gm$ for some $x_1\in \gm$. Note that
$x_1^p=0$. Changing $\d_1$ for $\d_1(x_1)^{-1}\d_1\in \Der (R)$,
one can assume that $\d_1(x_1)=1$.  Then, by Lemma
\ref{dpx27Jan06}.(3), one can assume that $\d_1^p=0$ (after
possibly changing $\d_1$). By Theorem \ref{30Jan06},
$R=\bigoplus_{i=0}^{p-1}R^{\d_1}x_1^i$. Note that $R^{\d_1}$ is a
local ring with maximal ideal $R^{\d_1}\cap \gm$. Let us prove
that {\em if $\ga $ is a differential ideal of the ring
$R^{\d_1}$, then $\ga':= \bigoplus_{i=0}^{p-1} \ga x_1^i$ is a
differential ideal of $R$}: for any $\eta \in \Der (R)$ and $c\in
R^{\d_1}$ we have $\eta (c) =\sum_{i=0}^{p-1} \eta_i( c) x_1^i$
where $\eta_i\in \Der (R^{\d_1})$, and the result follows. It
implies that $R^{\d_1}$ is a differentiably simple Noetherian
ring. Applying the same argument several times (or use Proposition
\ref{16Feb06} and induction) we have $R=\bigoplus_{0\leq i_\nu
\leq p-1} R^{\d_1, \ldots , \d_s} x_1^{i_1}\cdots x_s^{i_s}$ for
some {\em commuting} derivations $\d_1, \ldots , \d_s\in \Der (R)$
and elements $x_1, \ldots , x_s\in \gm$ such that $\d_i(x_j)
=\d_{ij}$ and $\d_1^p=\cdots =\d_s^p=0$, where $R^{\d_1, \ldots ,
\d_s}:=\bigcap_{i=1}^{s-1}\ker \, \d_i$ (note that any derivation
$\d $ of $R^{\d_1, \ldots , \d_s}$ can be  extended to a
derivation of $R$ by setting $\d (x_1)=\cdots = \d(x_s)=0$). Let
$s$ be the largest number for which there exist derivations $\d_1,
\ldots , \d_s$ as above (the $s$ exists since $\dim_k( \gm /
\gm^2)<\infty$). Then necessarily $R^{\d_1, \ldots , \d_s}$ is a
{\em field} canonically isomorphic to $k$,
 and $R=R^{\d_1, \ldots , \d_s}+\gm$.   Clearly, the set of elements
$x_1+\gm^2 ,\ldots , x_s+\gm^2 $ is a $k$-basis for $\gm / \gm^2$;
hence $n=s$.

$(3\Rightarrow 2)$: The determinant $\D := \det (\d_i(x_j))$ is a
unit of $R$. For each $i=1, \ldots , n$,  let us ``drop'' $x_i$ in
the determinant $\D $ and then multiply it by $\D^{-1}$; as the
result  we have well defined derivations of the ring $R$:
$$
\der_i (\cdot ) := \D^{-1} \det
 \begin{pmatrix}
  \d_1 (x_1) & \cdots & \d_n(x_1) \\
  \vdots & \vdots & \vdots \\
\d_1(\cdot ) & \cdots & \d_n(\cdot )\\
 \vdots & \vdots & \vdots \\
\d_1(x_n) & \cdots & \d_n(x_n) \\
\end{pmatrix},
\;\;\; i=1, \ldots , n,$$ such that $\der_i  (x_j)=\d_{ij}$ for
all $i,j=1,\ldots , n$.  Now, we finish the proof by applying
Proposition \ref{16Feb06}.

$7({\rm some})\Rightarrow 7({\rm any})$: Clearly, statement
$7$(some) implies statement $1$, and, as we have proved, statement
$1$ implies statement $2$. We can assume that $R= k[x_1, \ldots ,
x_n]/(x_1^p, \ldots , x_n^p)$ where $\gm = (x_1, \ldots , x_n)$.
Let $l$ be a coefficient field of $R$. Then $R= l+\gm = k+\gm$ and
$l\simeq R/\gm \simeq k$. Hence, for each $i\geq 0$,
$\dim_k(\gm^i/ \gm^{i+1})=\dim_{R/\gm}(\gm^i/
\gm^{i+1})=\dim_l(\gm^i/ \gm^{i+1})$. Therefore, $R\simeq l[x_1,
\ldots , x_n]/(x_1^p, \ldots , x_n^p)$. Now, it is obvious that
$R$ is a $\Der_l(R)$-simple $l$-algebra by using the partial
$l$-derivatives. $\Box $

 The next result gives {\em
explicitly} a subfield $k'$ of $R$ such that $R=k'+\gm$.

\begin{corollary}\label{cHB10Feb06}%\marginpar{cHB10Feb06}
Let $R$ be a differentiably simple Noetherian $\Fp$-algebra, the
derivations  $\d_1, \ldots , \d_n$ and the elements $x_1, \ldots ,
x_n\in \gm $ be as in Theorem \ref{Harpbig}.(4) (i.e.,  $\d_i^p=0$
and $\d_i(x_j)=\d_{ij}$ for all $i,j$). For each $i=1, \ldots ,
n$, let  $\phi_i:= \sum_{k=0}^{p-1}(-1)^k\frac{x_i^k}{k!}\d_i^k:
R\ra R$ and $\CN_n:= \{ \alpha = (\alpha_1, \ldots , \alpha_n)\in
\mathbb{N}^n\, | \, 0\leq \alpha_i \leq p-1$ for all $i\}$.
 Then $k':= {\rm im} \,
\phi$ is a coefficient  field for $R$ where $\phi$ is equal to the
composition of maps  $\prod_{i=1}^n\phi_i=\sum_{\alpha \in \CN_n}
(-1)^\alpha x^{[\alpha ]}\d^\alpha$ where  $x^{[\alpha ]}:=
\prod_{i=1}^n\frac{x_i^{\alpha_i}}{\alpha_i!}$, $\d^\alpha :=
\d_1^{\alpha_1}\cdots \d_n^{\alpha_n}$ and $(-1)^{\alpha }:=
(-1)^{\alpha_1 }\cdots (-1)^{\alpha_n }$.
\end{corollary}

{\it Proof}. We have proved already that $R=\bigoplus_{\alpha \in
\CN_n}k'x^{[\alpha ]}$ (see the proof $(1\Rightarrow 2)$  of
Theorem \ref{Harpbig}) where $k':= \bigcap_{i=1}^n \ker \, \d_i$
is a subfield of $R$ such that $R=k'+\gm$. By Theorem
\ref{30Jan06}, the map $\phi$ is a projection onto $k'$. $\Box $

{\it Remark}. In  view of Theorem \ref{Harpbig}.(3)-(6), Corollary
\ref{cHB10Feb06}, in fact, gives all the coefficient fields for
$R$ when combined with Proposition \ref{16Feb06} provided one
knows explicitly generators for the $R$-module $\Der (R)$. In more
detail,  one can find derivations $\d_1, \ldots , \d_n$ satisfying
Theorem \ref{Harpbig}.(3); then applying Proposition \ref{16Feb06}
one obtains derivations satisfying Corollary \ref{cHB10Feb06}
which produce the coefficient field $k'$. Varying the derivations
$\d_1, \ldots , \d_n$  one obtains all the coefficient fields for
$R$ as follows from Theorem \ref{b10Feb06}.

%%%%%%%%%%%%%%%%%% SECTION 3 %%%%%%%%%%%%%%%%%%%%%%%%

\section{Existence and uniqueness of an iterative  $\d$-descent}

The main result of this section is Theorem \ref{f18Jan06} on {\em
existence} and {\em uniqueness} of an {\em iterative} $\d$-{\em
descent}. This is the key  (and the most difficult) result of the
paper.

Let $\d $ be a derivation of a ring $R$. A finite sequence of
elements in $R$, $\yy $: $y^{[-1]}:=0, y^{[0]}:=1, y^{[1]}, \ldots
, y^{[m]}$ ($m\geq 1$) is called a $\d$-{\bf descent} if $\d
(y^{[i]})=y^{[i-1]}$ for all $i\geq 0$. If $R$ is an
$\Fp$-algebra, $m=p^n-1$, then a sequence in $R$,  $\{ y^{[i]},
0\leq i <p^n\}$,  is called an {\em
 iterative sequence}, if   $$y^{[i]}y^{[j]}={i+j\choose
i}y^{[i+j]}, \;\;\; {\rm  for\; all}\;\;  0\leq i,j\leq p^n-1.$$
Note that $y^{[i]}y^{[j]}=0$ if $i+j\geq p^n$ since ${i+j\choose
i}=0$ in $\Fp$.

{\it Definition}. An iterative sequence $\{ y^{[i]}, 0\leq i
<p^n\}$ which is a $\d$-descent is called an {\bf iterative}
$\d$-{\bf descent} of {\em exponent} $n$. Note that any truncation
 $\{ y^{[i]}, 0\leq i <p^m\}$, $1\leq m <n$, of the iterative
 $\d$-descent $\{ y^{[i]}, 0\leq i
<p^n\}$ is  an iterative $\d$-descent of exponent $m$.

The following lemma establishes relations between iterative
descents and simple derivations.

\begin{lemma}\label{Ky13Feb06}%\marginpar{Ky13Feb06}
Let $\d$ be a derivation of an $\Fp$-algebra $R$, $K:= \ker \,
\d$, and $\xx =\{ x^{[i]}, 0\leq i <p^n\}$ be an iterative
$\d$-descent, elements of which commute with $K$.  Then
\begin{enumerate}
\item the $K$-algebra $K\langle \xx \rangle$ generated over $K$ by
all the elements $x^{[i]}$ is equal to
$$K\langle x_0, \ldots ,
x_{n-1}\rangle=\bigoplus_{i=0}^{p^n-1}Kx^{[i]}\simeq K[x_0, \ldots
, x_{n-1}]/(x_0^p,\ldots ,x_{n-1}^p)$$ where $x_k:= x^{[p^k]}$ and
$x^{[i]}= \prod_{k=0}^t \frac{x_k^{i_k}}{i_k!}$ where
$i=\sum_{k=0}^ti_kp^k$, $0\leq i_k<p$, the $p$-adic presentation
of the integer $i$.
 \item $\d':=\d|_{K\langle \xx \rangle}=\sum_{k=0}^{n-1}x^{[p^k-1]}\frac{\der }{\der x_k}\in \Der_K(K\langle \xx \rangle )$.
 \item If $K$ is a field, then $\d'$ is a simple $K$-derivation of
 the algebra $K\langle \xx \rangle$.
\end{enumerate}
\end{lemma}

{\it Proof}. $1$. Clearly, the equality  $K\langle \xx \rangle
=\sum_{i=0}^{p^n-1}Kx^{[i]}$ holds  since $\xx$ is an iterative
sequence. Since the sequence $\xx$ is a $\d$-descent, it follows
easily that the sum is a direct one, i.e.,  $K\langle \xx \rangle
=\bigoplus_{i=0}^{p^n-1}Kx^{[i]}$. The $\xx$ is an iterative
sequence, hence $x_0^p=\cdots = x_{n-1}^p=0$ and $x^{[i]}=
\prod_{k=0}^t \frac{x_k^{i_k}}{i_k!}$ where
$i=\sum_{k=0}^ti_kp^k$, $0\leq i_k<p$. So, there is a natural
$K$-algebra isomorphism $K[x_0, \ldots , x_{n-1}]/(x_0^p,\ldots
 , x_{n-1}^p)\simeq K\langle \xx \rangle$.

$2$. This is obvious.

$3$. If $a= a_0+a_1x^{[1]}+\cdots +a_sx^{[s]}$ is a nonzero
element of the algebra $K\langle \xx \rangle$ where $a_i\in K$ and
$a_s\neq 0$, then $\d^s(a_s^{-1}a)=1$. Therefore, $\d'$ is a
simple $K$-derivation of the algebra $K\langle \xx \rangle$.
$\Box $

An algebra $S$ over a field $K$ is a positively filtered algebra
if $S$ is a union of its subspaces, $S= \bigcup_{i\geq 0} S_i$,
such that $K\subseteq    S_0\subseteq S_1\subseteq \cdots $ and
$S_iS_j \subseteq S_{i+j}$ for all $i,j\geq 0$.  Let $A$ be an
algebra over a field $K$ and let $\d $ be a
 $K$-derivation of the algebra $A$. For any elements $a,b\in A$
 and a natural number $n$, an easy induction argument yields
 $$ \d^n(ab)=\sum_{i=0}^n\, {n\choose i}\d^i(a)\d^{n-i}(b).$$
 It follows that the union of the
 vector spaces $N:=N(\d ,A)=\bigcup_{i\geq 0}\, N_i$, $N_i:= \ker \, \d^{i+1}$,  is a positively
 {\em filtered} algebra ($N_iN_j\subseteq N_{i+j}$ for
 all $i,j\geq 0$), so-called, the {\em nil algebra} of $\d $.
 Clearly,  $N_0= A^\d:=\ker \, \d $  is a
 subalgebra (of {\em constants} for $\d $) of $A$,
   and $N=\{ a\in A \, | \ \d^n (a)=0$
  for some natural $n=n(a)\}$.

\begin{lemma}\label{fyi18Jan06}%\marginpar{fyi18Jan06}
Let  $\d $ be a derivation of a ring $A$ and $\{ x^{[i]}, 0\leq
i\leq m\}$ be a $\d$-descent. Then $N(\d ,
A)_i=\bigoplus_{j=0}^iA^\d x^{[j]}= \bigoplus_{j= 0}^ix^{[j]}
A^\d$, $0\leq i\leq m$.
\end{lemma}

{\it Proof}. For each $i\geq 0$, let $N_i':=\bigoplus_{j=0}^i A^\d
x^{[j]}$ and $ N_i:= N(\d , A)_i$. Then $N_i'\subseteq N_i$,
$i\geq 0$.  Clearly, $N_0'=A^\d = N_0$. We use induction on $i$ to
prove that
 $N_i'=N_i$. Let $i\geq 1$ and $N_{i-1}'=N_{i-1}$ (by the
 induction hypothesis). Let $u\in N_i$. Then $c:= \d^i(u)\in A^\d$; hence
$\d^i(u) = \d^i (cx^{[i]})$, and so $u-cx^{[i]}\in N_{i-1}=
N_{i-1}'$. Therefore, $u\in N_i'$, and so $N_i'=N_i$. Since $\d $
is a derivation of the opposite
 algebra $A^{op}$,  the equalities $N_i=\bigoplus_{j=0}^i x^{[j]}A^\d
$, $i\geq 0$, follow from the just proved  equalities. $\Box $

\begin{lemma}\label{ft2Dec05}%\marginpar{ft2Dec05}
Let  $\d $ be a derivation of a ring $A$ and $\{ x^{[i]}, 0\leq
i\leq m\}$ be a $\d$-descent.  Then $\{ x^{[i]'}, 0\leq i\leq m\}$
is a $\d$-descent iff $x^{[0]'}:=1$ and
$x^{[i]'}=x^{[i]}+\sum_{j=1}^i\l_jx^{[i-j]}$, $1\leq i\leq m$,
where all $\l_j\in A^\d$.
\end{lemma}

{\it Proof}. $(\Leftarrow )$ Obvious.

$(\Rightarrow )$ We prove this implication by induction on $i\geq
1$. Let $i=1$. Then $\d (x^{[1]'})=1=\d (x^{[1]})$ implies
$x^{[1]'}=x^{[1]}+\l_1$ for some element $\l_1\in A^\d$. Suppose
that $i\geq 2$, and, by the induction hypothesis,
$x^{[i-1]'}=x^{[i-1]}+\sum_{j=1}^{i-1}\l_jx^{[i-1-j]}$ for some
$\l_j\in A^\d$. Then $\d (x^{[i]'})=x^{[i-1]'}=\d
(x^{[i]}+\sum_{j=1}^{i-1}\l_jx^{[i-j]})$ implies $\l_i:=
x^{[i]'}-x^{[i]}-\sum_{j=1}^{i-1}\l_jx^{[i-j]}\in A^\d$, as
required.  $\Box $

\begin{lemma}\label{fcxi18Jan06}%\marginpar{fcxi18Jan06}
Let $\d$ be a derivation of a ring  $A$ such that $\d^i
(y^{[i]})=1$, $0\leq i\leq m$, for some elements $y^{[i]}$ of $A$.
 Note that $y^{[0]}=1$. Then there exists a {\bf unique}
$\d$-descent $ \{ x^{[i]},  0 \leq i\leq m\}$  such that
 $x^{[1]}=y^{[1]}$ and
 $x^{[i]}=y^{[i]}+\sum_{j=1}^{i-1}c_{ij}y^{[j]}$, $2\leq i\leq m $, for some  $c_{ij}\in A^\d$.
\end{lemma}

{\it Proof}. One can easily prove that $N_i:= N(\d ,
A)_i=\bigoplus_{j=0}^iA^\d y^{[j]}$, $0\leq i \leq m$ (repeat the
argument of the proof of Lemma \ref{fyi18Jan06}). Let, for a
moment, a sequence  $ \{ x^{[i]'}, 0 \leq i\leq m\}$ be an
arbitrary  $\d$-descent. Then, by Lemma \ref{fyi18Jan06},
$x^{[i]'}\in N_i$, i.e.,
$x^{[i]'}=y^{[i]}+\sum_{j=0}^{i-1}c_{ij}'y^{[j]}$, $1\leq i\leq
m$, for some elements $c_{ij}'\in A^\d$ (note that one can easily
find a $\d$-descent, e.g.  $\{ z^{[i]}:= \d^{m-i}(y^{[m]}), 0\leq
i \leq m\})$. Let $ \{ x^{[i]}, 0 \leq i\leq m\}$ be another
$\d$-descent, and so
$x^{[i]}=y^{[i]}+\sum_{j=0}^{i-1}c_{ij}y^{[j]}$, $1\leq i\leq m$,
for some elements $c_{ij}\in A^\d$. By Lemma \ref{ft2Dec05},
$$x^{[i]}=x^{[i]'}+\sum_{j=1}^i\l_jx^{[i-j]'}, \;\; 1\leq i\leq
m,$$ for some elements $\l_j\in A^\d$.
 We have to prove that the defining conditions
$$c_{1,0}=c_{2,0}=\cdots =c_{m,0}=0$$ of the $\d$-descent from
Lemma \ref{fcxi18Jan06} {\em uniquely} determine the elements
$\l_1, \l_2, \ldots , \l_m$. The equality $c_{1,0}=0$ yields the
equalities  $y^{[1]}=x^{[1]}=x^{[1]'}+\l_1=y^{[1]}+c_{1,0}'+\l_1$;
hence $\l_1=-c_{1,0}'$. Suppose that, using the equalities
$c_{1,0}=\cdots =c_{i-1,0}=0$, we have already found unique
elements $\l_1,\l_2,\ldots , \l_{i-1}$. Then the element $\l_i$
can be found uniquely from the equality
$x^{[i]}=x^{[i]'}+\l_1x^{[i-1]'}+\cdots +\l_{i-1}x^{[1]'}+\l_i$.
We have to equate to zero the coefficient $c_{i,0}$  of
$y^{[0]}:=1$ after we substitute the sum for each $x^{[i]'}$ above
(via $y^{[k]}$):
$$
x^{[i]}=y^{[i]}+\sum_{j=1}^{i-1}c_{ij}y^{[j]}+c_{i,0}'+\l_1c_{i-1,0}'+\cdots
+\l_{i-1}c_{1,0}'+\l_i;$$ that is, $\l_i:=
-c_{i,0}'-\l_1c_{i-1,0}'-\cdots -\l_{i-1}c_{1,0}'$. Therefore,
 for this unique choice of $\{ \l_i\}$, we have $c_{1,0}=c_{2,0}=\cdots
 =c_{m,0}=0$ for the $\d$-descent $\{ x^{[i]}\}$ in Lemma \ref{fcxi18Jan06}. $\Box$

{\bf Binomial coefficients modulo $p$}. For any two nonnegative
integers $i$ and $j$ written in the $p$-adic form as
$i=\sum_{k}i_kp^k$, $0\leq i_k<p$,  and $j=\sum_kj_kp^k$, $0\leq
j_k<p$,  in the
field $\Fp$ there is the equality %\marginpar{ijbin}
\begin{equation}\label{ijbin}
 {i\choose
j}=\prod_k{i_k\choose j_k}.
\end{equation}
The equality is obvious if $i<j$ since both sides of the equality
are equal to zero. If $i\geq j$, this equality can be proved by
comparing the coefficients of
 $x^j$ of the polynomials in $\Fp [x]$ at both ends of the equality
 $$ \sum_{j=0}^i {i\choose j} x^j= (1+x)^i =
 \prod_k(1+x^{p^k})^{i_k}= \prod_k (\sum_{j_k=0}^{i_k} {i_k\choose
 j_k} x^{j_kp^k}) = \sum_{j=0}^i \prod_k {i_k\choose j_k} x^j.$$

It follows that in $\Fp $,  %\marginpar{1ijbin}
\begin{equation}\label{1ijbin}
{i\choose j }\neq 0 \;\; {\rm iff}\;\; j_k\leq i_k\;\; {\rm for \;
all}\;k,
\end{equation}
%\marginpar{2ijbin}
\begin{equation}\label{2ijbin}
 {i+j\choose j}\neq 0\;\; {\rm  iff}\;\;  i_k+j_k<p\;\; {\rm for \;
 all}\;\; k, \;\; {\rm and}
\end{equation}
%\marginpar{3ijbin}
\begin{equation}\label{3ijbin}
{ip^s\choose jp^s }={i\choose j}, \;\; s\geq 1.
\end{equation}

Let $A$ be an $\Fp$-algebra. Recall that a  sequence  $\{x^{[i]},
0\leq i <p^n\}$ in $A$ is called  an {\em iterative sequence} iff
$x^{[i]}x^{[j]}= {i+j\choose j}x^{[i+j]}$ for all $0\leq i,j<p^n$
where $x^{[k]}:=0$ for $k\geq p^n$. Note that if $i+j\geq p^n$,
then $x^{[i]}x^{[j]}= {i+j\choose j}x^{[i+j]}=0\cdot x^{[i+j]}=0.
$

\begin{proposition}\label{ittp}%\marginpar{ittp}
{\rm (Structure of iterative sequence)}  Let $A$ be an
$\Fp$-algebra and $\{x^{[i]}, 0\leq i <p^n\}$ be an iterative
sequence. Then
\begin{enumerate}
\item for each $i=1, \ldots , p^n-1$, written $p$-adically as $i=
\sum_k i_kp^k$, $x^{[i]}= \prod_k \frac{x^{[p^k]i_k}}{i_k!}$. This
means that the iterative sequence is determined by the elements
$\{ x^{[0]}, x^{[p^j]} \, | \,  j=0,1, \ldots ,  n-1\}$.  \item
For each $j=0,1,   \ldots, n-1$, $x^{[p^j]p}=0$ (hence
$x^{[i]p}=0$ for all $i=1, \ldots , p^n-1$, by statement 1). \item
$x^{[0]}x^{[p^j]}= x^{[p^j]}$, $j=0, 1, \ldots , n-1$, and
$x^{[0]}x^{[0]}= x^{[0]}$.
\end{enumerate}
Conversely, given commuting elements $\{ x^{[0]}, x^{[p^j]}\, | \,
j=0,1, \ldots ,  n-1\}$, in $A$ that satisfy the conditions of
statements 2 and 3 above,  then the elements $\{ x^{[i]}, 0\leq i
<p^n\}$ defined as in statement 1 form an iterative sequence.
\end{proposition}

{\it Remark}. To make formulae more readable we often use the
notation $x^{[p^k]j}$ for $(x^{[p^k]})^j$.

{\it Proof}. 1. We have to prove that $\prod_k
\frac{x^{[p^k]i_k}}{i_k!}= x^{[\sum_k i_kp^k]}$. Consider first a
special case using (\ref{3ijbin}),
$$ \frac{x^{[p^k]i_k}}{i_k!}= \frac{{2p^k\choose p^k} {3p^k\choose p^k}\cdots {i_kp^k\choose p^k}}{i_k!} x^{[i_kp^k]}
=\frac{{2\choose 1} {3\choose 1}\cdots {i_k\choose 1}}{i_k!}
x^{[i_kp^k]}= \frac{i_k!}{i_k!} x^{[i_kp^k]}=x^{[i_kp^k]}.
$$
Now, the general case follows from the special case and
(\ref{ijbin}) by simply multiplying the elements below and using
the fact that  each multiplication yields a binomial which is 1 in
$\Fp$, by (\ref{ijbin}):
$$\prod_{k=0}^s \frac{x^{[p^k]i_k}}{i_k!}=\prod_{k=0}^s x^{[i_kp^k]}=
 \cdots =
x^{[\sum_{k=0}^{t-1}i_kp^k]}x^{[i_tp^t]}\cdots x^{[i_sp^s]}=\cdots
= x^{[\sum_{k=0}^si_kp^k]}=x^{[i]}.$$

2. For each $i=1, \ldots , p^n-1$,
$$ (x^{[i]})^p= x^{[i]}x^{[i]}\cdots x^{[i]}= {2i\choose
i}{3i\choose i}\cdots {pi\choose i}x^{[pi]}=0\cdot x^{[pi]}=0,$$
since ${pi\choose i}=0$ in $\Fp$.

3. This is obvious.

Conversely, suppose that elements $\{ x^{[0]}, x^{[p^j]}\, | \,
j=0,1, \ldots , n-1\}$  satisfy the conditions of statements 2 and
3, and that the elements $\{ x^{[i]}, 0\leq i <p^n\}$ are defined
as in statement 1. To prove that the sequence $\{ x^{[i]}, 0\leq i
<p^n\}$ is iterative it suffices to show that
$x^{[i]}x^{[j]}={i+j\choose j}x^{[i+j]}$ for all $1\leq i,j<p^n$.
Let $i=\sum i_kp^k$ and $j=\sum j_kp^k$ be the $p$-adic forms of
$i$ and $j$. Suppose that $i_k+j_k\geq p$ for some $k$. Then, on
the one hand, ${i+j\choose j}=0$ (by (\ref{2ijbin})), and so
$x^{[i]}x^{[j]}=0={i+j\choose j}x^{[i+j]}$. Suppose that
$i_k+j_k<p$ for all $k$. Then
\begin{eqnarray*}
x^{[i]}x^{[j]}&=&\prod_k\frac{(x^{[p^k]})^{i_k}}{i_k!}\frac{(x^{[p^k]})^{j_k}}{j_k!}=\prod_k{i_k+j_k\choose
i_k}\frac{(x^{[p^k]})^{i_k+j_k}}{(i_k+j_k)!} \\
&=& \prod_k{i_k+j_k\choose i_k}\cdot \prod_l
\frac{(x^{[p^l]})^{i_l+j_l}}{(i_l+j_l)!}= {i+j\choose i}x^{[i+j]}.
\end{eqnarray*}
This means that $\{ x^{[i]}, 0\leq i <p^n\}$ is an iterative
sequence. $\Box $

The following corollary gives  necessary and sufficient conditions
for an iterative sequence to be a $\d$-descent.

\begin{corollary}\label{1ittp}%\marginpar{1ittp}
Let $A$ be an $\Fp$-algebra, $\d$ be a derivation of $A$, and $\{
x^{[i]}, 0\leq i<p^n\}$ be an iterative sequence in $A$ with
$x^{[0]}=1$. Then the iterative sequence $\{ x^{[i]}, 0\leq
i<p^n\}$ is a $\d$-descent iff $\d (x^{[p^j]})= x^{[p^j-1]}$,
$0\leq j \leq n-1$.
\end{corollary}

{\it Proof}. $(\Rightarrow )$ Trivial.

$(\Leftarrow )$ We have to show that $\d (x^{[i]})= x^{[i-1]}$,
$1\leq i<p^n$. Observe that, for each $k\geq 1$,
$p^k-1=\sum_{l=0}^{k-1}(p-1)p^l$. Then, by (\ref{2ijbin}),
${p^j-1+p^s\choose p^s}=0$, $0\leq s<j$, and so
$$ x^{[p^s]}\cdot \d (x^{[p^j]})= x^{[p^s]}x^{[p^j-1]}={p^j-1+p^s\choose
p^s}x^{[p^s+p^j-1]}=0.$$ These equalities imply that, for any
integer $i$ $(1\leq i<p^n)$ written $p$-adically as
$i=i_sp^s+i_{s+1}p^{s+1}+\cdots +i_tp^t$ with $i_s\neq 0$, $s\leq
t$,  and $i=i_sp^s+j$, $j:= i_{s+1}p^{s+1}+\cdots +i_tp^t$,
\begin{eqnarray*}
 \d (x^{[i]})& =& \d (x^{[i_sp^s]}x^{[j]})= \d (\frac{ x^{[p^s]i_s}  }{i_s!} x^{[j]})=\d (\frac{ x^{[p^s]i_s}  }{i_s!} ) x^{[j]}=
 x^{[p^s-1]}\frac{ x^{[p^s](i_s-1)}  }{(i_s-1)!}x^{[j]}\\
 &=&x^{[p^s-1]} x^{[p^s(i_s-1)]}x^{[j]}=x^{[p^si_s-1]}
 x^{[j]}=x^{[p^si_s-1+j]}=x^{[i-1]}.\;\; \Box
\end{eqnarray*}

Combining Proposition \ref{ittp} and Corollary \ref{1ittp}, one
obtains necessary and sufficient conditions for a sequence to be
an iterative $\d$-descent.

\begin{corollary}\label{2ittp}%\marginpar{2ittp}
Let $A$ be an $\Fp$-algebra, $\d$ be a derivation of $A$, and
$x^{[1]}, x^{[p]},  \ldots , x^{[p^{n-1}]}$ be commuting elements
of $A$. Let  $x^{[i]}:=\prod_k  \frac{x^{[p^k]i_k}}{i_k!}$ for
each $i=\sum_k i_kp^k$, $0\leq i_k<p$, such that $0\leq i<p^n$.
Then the sequence $ \{ x^{[i]}, 0\leq i <p^n\}$ is an iterative
$\d$-descent iff $\d ( x^{[p^j]})= x^{[p^j-1]}$ and $
x^{[p^j]p}=0$ for $0\leq j \leq n-1$.
\end{corollary}

{\it Proof}. $(\Rightarrow )$ Trivial.

$(\Leftarrow )$ The conditions $ x^{[p^j]p}=0$, $0\leq j \leq
n-1$, mean that $ \{ x^{[i]}, 0\leq i <p^n\}$ is an iterative
sequence, by Proposition \ref{ittp}. Then, the conditions $\d (
x^{[p^j]})= x^{[p^j-1]}$, $0\leq j \leq n-1$, imply that $ \{
x^{[i]}, 0\leq i <p^n\}$ is a $\d$-descent, by Corollary
\ref{1ittp}. $\Box$

Let $A$ be a commutative $\Fp$-algebra and $\d$ be a derivation of
$A$. Let $\ID (\d, n)$ be the set of all iterative $\d$-descents $
\{ x^{[i]}, 0\leq i <p^n\}$ of exponent $n$ in $A$. Let $C(\d ,
n)$ be the set of all $n$-tuples $(\l_0, \l_1, \ldots , \l_{n-1})$
such that $\l_j\in A^\d$ and $\l_j^p=0$ for  $0\leq j\leq n-1$.
Note that if $A^\d$ is a {\em reduced} ring, then $C(\d , n) = \{
(0, \ldots , 0)\}$, i.e.,  $C(\d , n)$ contains a single element.
 By Lemma \ref{fyi18Jan06} and Proposition \ref{ittp}, for each iterative $\d$-descent,
 say $ \{ x^{[i]}, 0\leq i
<p^n\}$,
%\marginpar{NdApn}
\begin{equation}\label{NdApn}
N(\d, A)_{p^n-1}= \bigoplus_{i=0}^{p^n-1}A^\d x^{[i]}\simeq A^\d [
x^{[1]}, x^{[p]}, \ldots , x^{[p^{n-1}]}]/(x^{[1]p}, x^{[p]p},
\ldots , x^{[p^{n-1}]p}).
\end{equation}
 So, $N(\d,
A)_{p^n-1}$ is a {\em subring} of $A$ that contains $A^\d$, and
the decomposition (\ref{NdApn}) holds for {\em all} iterative
$\d$-descents in $A$ of exponent $n$. In particular, all iterative
$\d$-descents in $A$ of exponent $n$ belong to $N(\d, A)_{p^n-1}$.
 Then
%\marginpar{NNdA}
\begin{equation}\label{NNdA}
N(\d, A)_{p^n-1}=A^\d \oplus \gm, \;\; \gm :=(x^{[1]}, x^{[p]},
\ldots , x^{[p^{n-1}]}).
\end{equation}
If $\{ y^{[i]}, 0\leq i <p^n\}$ is an iterative $\d$-descent in
$A$, then $\{ y^{[i]}, 0\leq i <p^n\}\subseteq N(\d, A)_{p^n-1}$.
Therefore, the following map is well defined: %\marginpar{rIDC}
\begin{equation}\label{rIDC}
r=r_n: \ID (\d , n)\ra C(\d , n), \;\; \{ y^{[i]}, 0\leq i
<p^n\}\mapsto (\l_0, \l_1, \ldots , \l_{n-1}),
\end{equation}
where $\l_j\equiv y^{[p^j]}\mod \gm$, $j=0, 1, \ldots , n-1$. Note
that the map $r$ depends on the choice of the iterative
$\d$-descent $\{ x^{[i]}, 0\leq i <p^n\}$ since the decomposition
(\ref{NNdA}) does.

\begin{theorem}\label{f18Jan06}%\marginpar{f18Jan06}
{\rm (Existence and uniqueness of an iterative $\d$-descent)}  Let
$A$ be a commutative algebra over a field $K$ of characteristic
$p>0$ and $\d $ be a $K$-derivation of the algebra $A$ such that
there exists a finite sequence of elements $y_0, y_1, \ldots ,
y_{n-1}$ of $A$ such that $y_k^p=0$ and $\d^{p^k}(y_k)=1$ for all
$0\leq k\leq n-1$. Then
\begin{enumerate}
\item (Existence) The following sequence $\{ x^{[i]}, 0\leq i
<p^n\}$ is an iterative $\d$-descent where $x^{[0]}:=1$,
$x^{[1]}:=y_0$,  and, for  $i\geq 2$ written $p$-adically as
$i=\sum_{k=0}^ti_kp^k$ ($0\leq i_k\leq p-1$) the element $x^{[i]}$
is defined as
$x^{[i]}:=\prod_{k=0}^t\frac{(x^{[p^k]})^{i_k}}{i_k!}$, where
$$x^{[p]}:= (-1)^{p-1}\phi_0(y_1), \;\; \phi_0 (z):=
\sum_{j=0}^{p-1}(-1)^j\frac{(x^{[1]})^j}{j!}\d^j(z),$$ and then
recursively, for each $k$ such that  $1\leq k \leq n-2$, the
element $x^{[p^{k+1}]}$ is defined by the rule
$$x^{[p^{k+1}]}:= (-1)^{p-1}\d^{p^k-1} (\prod_{l=0}^{k-1}
\frac{(x^{[p^l]})^{p-1}}{(p-1)!}\cdot \phi_k(y_{k+1})), \;\;
\phi_k (z):=
\sum_{j=0}^{p-1}(-1)^j\frac{(x^{[p^k]})^j}{j!}\d^{p^kj}(z).$$
\item (Almost uniqueness) Let $\{ x^{[i]}, 0\leq i <p^n\}$ be an
arbitrary  iterative $\d$-descent (not necessarily as in statement
1, and $n$ here is not necessarily as in statement 1 either). Then
the map (\ref{rIDC}) is a bijection. \item (Uniqueness). If, in
addition, the ring $A^\d$ is reduced, then $\{ x^{[i]}, 0\leq i
<p^n\}$ from statement 1  is the only iterative $\d$-descent.
\end{enumerate}
\end{theorem}

 {\it Proof}. 1. By Corollary \ref{2ittp},
 it suffices to prove two statements: $\d (x^{[p^j]})=x^{[p^j-1]}$ and
$(x^{[p^j]})^p=0$ for all $0\leq j\leq n-1$. For  $j=0$, we have
$\d (x^{[1]})=\d (y_0)=1=x^{[0]}$ and $(x^{[1]})^p=y_0^p=0$. A
direct calculation shows that $\d \phi_0
(z)=(-1)^{p-1}\frac{y_0^{p-1}}{(p-1)!}\d^p(z)$. If $j=1$, then $\d
(x^{[p]})=(-1)^{p-1}\d
\phi_0(y_1)=(-1)^{p-1}(-1)^{p-1}\frac{y_0^{p-1}}{(p-1)!}\d^p(y_1)=x^{[p-1]}$
and
$$
(x^{[p]})^p=(-1)^{(p-1)p}(\sum_{k=0}^{p-1}(-1)^k\frac{y_0^k}{k!}\d^k(y_1))^p=
\sum_{k=0}^{p-1}(-1)^{kp}(\frac{y_0^k}{k!})^p(\d^k(y_1))^p=y_1^p=0.$$

Let $j\geq 2$. We use induction on $j$.  By the induction
hypothesis, for all $k<j$, $\d (x^{[p^k]})=x^{[p^k-1]}$ and
$(x^{[p^k]})^p=0$. This means that $\{ x^{[i]}, 0\leq i <p^j\}$ is
an iterative $\d$-descent,  by Corollary \ref{2ittp}.  In
particular, $\d (x^{[l]})=x^{[l-1]}$ for all
$l<p^j$, which implies that %\marginpar{1RR}
\begin{equation}\label{1RR}
x^{[p^j-1]}=x^{[\sum_{l=0}^{j-1}(p-1)p^l]}=\prod_{l=0}^{j-1}\frac{(x^{[p^l]})^{p-1}}{(p-1)!}\in
\ker \, \d^{p^j}.
\end{equation}
 For each $k$ such that $2\leq k<j$, a direct calculation shows
that %\marginpar{2RR}
\begin{equation}\label{2RR}
\d^{p^k}\phi_k(z)=
(-1)^{p-1}\frac{(x^{[p^k]})^{p-1}}{(p-1)!}\d^{p^{k+1}}(z).
\end{equation}
Now, using the two equalities above we have
\begin{eqnarray*}
 \d (x^{[p^j]})&=&
 (-1)^{p-1}\d^{p^{j-1}}(\prod_{l=0}^{j-2}\frac{(x^{[p^l]})^{p-1}}{(p-1)!}\phi_{j-1}(y_j))\\
 &=& (-1)^{p-1}\d^{p^{j-1}}(x^{[p^{j-1}-1]}\phi_{j-1}(y_j))\\
 &=&(-1)^{p-1}x^{[p^{j-1}-1]}\d^{p^{j-1}}\phi_{j-1}(y_j)\;\;\;\;\;\;\;\;\; \;\;\;\;\;\;\;\;\;\;\;\;\;\;\;\;  ({\rm by} \; \ref{1RR}) \\
 &=&(-1)^{p-1}x^{[p^{j-1}-1]}(-1)^{p-1}\frac{(x^{[p^{j-1}]})^{p-1}}{(p-1)!}\d^{p^j}(y_j)\;\;\;\;\;\; ({\rm by} \; \ref{2RR})\\
 &=&\prod_{l=0}^{j-1}\frac{(x^{[p^l]})^{p-1}}{(p-1)!}=x^{[p^j-1]}.
\end{eqnarray*}
Finally, letting  $t:=p^{j-1}-1$,
\begin{eqnarray*}
 x^{[p^j]}&=&
 (-1)^{p-1}\d^t(x^{[t]}\phi_{j-1}(y_j))=(-1)^{p-1}\sum_{s=0}^t{t\choose
 s}\d^s(x^{[t]})\d^{t-s}\phi_{j-1}(y_j)\\
 &=& (-1)^{p-1}\sum_{s=0}^t{t\choose
 s}x^{[t-s]}\d^{t-s}\phi_{j-1}(y_j).
\end{eqnarray*}
Since $(x^{[1]})^p=\cdots = (x^{[t]})^p=0$, it follows from the
equality above that $(x^{[p^j]})^p=0$ iff $\phi_{j-1}(y_j)^p=0$.
Since
$\phi_{j-1}(y_j)=\sum_{k=0}^{p-1}(-1)^k\frac{(x^{[p^{j-1}]})^k}{k!}\d^{kp^{j-1}}(y_j)$
and $y_j^p=0=(x^{[p^{j-1}]})^p$, we have $\phi_{j-1}(y_j)^p=0$.
This proves that $(x^{[p^j]})^p=0$, as required.

2. In order to prove statement 2, we use induction on $n\geq 1$.
The case $n=1$ is almost obvious. If $\{ y^{[i]}, 0\leq i <p\}\in
\ID (\d , 1)$, then $\d (y^{[1]})=1= \d (x^{[1]})$, and so
$y^{[1]}= x^{[1]}+\l_0$ for some element $\l_0\in A^\d$
necessarily $\l_0^p=0$ since $y^{[1]p}=x^{[1]p}=0$; and
$y^{[1]}\equiv \l_0\mod \gm$. By Corollary \ref{2ittp}, the
sequence $\{ y^{[i]}, 0\leq i<p\}$ is uniquely determined by the
element $y^{[1]}$, and so the map $r$ is injective. It remains to
show that $r$ is surjective. For each element, say $\l_0\in A^\d$,
 such that $\l_0^p=0$, the element $y^{[1]}= x^{[1]}+\l_0$
 satisfies the following conditions: $y^{[1]}\equiv \l_0\mod \gm$,
 $\d (y^{[1]})=1$, and $y^{[1]p}= x^{[1]p}+\l_0^p=0$. Hence, by
 Corollary \ref{2ittp}, the element $y^{[1]}$ determines an
 iterative $\d$-descent. Therefore, $r$ is a surjection, as
 required.

Now, let $n\geq 2$. Suppose that the result is true for all
$n'<n$. First, let us prove that the map $r:=r_n$ is injective.
Let $y:= \{ y^{[i]}, 0\leq i <p^n\}$ and $z:=\{ z^{[i]}, 0\leq i
<p^n\}$ be two iterative $\d$-descents such that $r(y) = r(z)=
(\l_0, \ldots , \l_{n-1})$. We have to show that $y=z$. By
induction, $y^{[i]}= z^{[i]}$, $0\leq i <p^{n-1}$. It follows from
the equalities
$$ \d (y^{[p^{n-1}]})= y^{[p^{n-1}-1]}= z^{[p^{n-1}-1]}=\d
(z^{[p^{n-1}]})$$ that $y^{[p^{n-1}]}-z^{[p^{n-1}]}\in A^\d$.
Since $y^{[p^{n-1}]}\equiv \l_{n-1}\equiv z^{[p^{n-1}]}\mod \gm$
and $N(\d , A)_{p^n-1} = A^\d \oplus \gm$ (by (\ref{NNdA})), we
must have $y^{[p^{n-1}]}=z^{[p^{n-1}]}$. By Corollary \ref{2ittp},
$y=z$, i.e.,  $r$ is an injection.

It remains to show that $r$ is a surjection. Let $\l := (\l_0,
\ldots , \l_{n-1})\in C(\d , n)$. We have to show that there
exists an element $y:= \{ y^{[i]}, 0\leq i <p^n\} \in \ID (\d ,
n)$ such that $r(y) = \l$. By induction, there exists a unique
element $y':= \{ y^{[i]}, 0\leq i <p^{n-1}\} \in \ID (\d , n-1)$
such that $r_{n-1}(y') = (\l_0, \ldots , \l_{n-2})$. By Corollary
\ref{2ittp}, it suffices to find an element $y^{[p^{n-1}]}$ such
that $y^{[p^{n-1}]p}=0$, $\d (y^{[p^{n-1}]})= y^{[p^{n-1}-1]}$,
and $y^{[p^{n-1}]}\equiv \l_{n-1}\mod \gm$. By Lemma
\ref{fyi18Jan06},
$$ N(A, \d )_{p^{n-1}-1}= \bigoplus_{i=0}^{p^{n-1}-1}A^\d x^{[i]}=
\bigoplus_{i=0}^{p^{n-1}-1}A^\d y^{[i]}, \;\; N(A, \d
)_{p^{n-1}}=N(A, \d )_{p^{n-1}-1}\oplus A^\d x^{[p^{n-1}]}.$$
Since the map $\d : N(A, \d )_{p^{n-1}}\ra N(A, \d )_{p^{n-1}-1}$
is {\em surjective} and $y^{[p^{n-1}-1]}\in N(A, \d
)_{p^{n-1}-1}$, one can find an element $y^{[p^{n-1}]}\in N(A, \d
)_{p^{n-1}}$ such that $\d (y^{[p^{n-1}]})= y^{[p^{n-1}-1]}$. The
element $y^{[p^{n-1}]}$ is unique up to adding an element of
$A^\d$. By adding a well chosen element of $A^\d$ to
$y^{[p^{n-1}]}$, we can assume that $y^{[p^{n-1}]}\equiv
\l_{n-1}\mod \gm$; i.e.,  $y^{[p^{n-1}]}=\l_{n-1}+v$ for some
element $v\in \gm$. Since $\l_{n-1}^p=0$ and $v^p=0$ (since $v\in
\gm )$, $y^{[p^{n-1}]p}=0$. Now, by Corollary \ref{2ittp}, the
elements $y^{[p^j]}$, $0\leq j \leq n-1$, determine an element,
say $y$, of $\ID (A, n)$ such that, obviously, $r(y)=\l$. This
proves that $r$ is a surjection. By induction, statement 2 holds.

3. Since $A^\d$ is a reduced ring, the set $C(\d ,n)$ contains the
single element $(0, \ldots , 0)$, and so the result follows from
statement 2. $\Box$

The next lemma shows that the set $\nsder (R)$ is nonempty where
$R$ is a differentiably simple Noetherian commutative ring.

\begin{lemma}\label{NnsdR}%\marginpar{NnsdR}
Let $k$ be a field of characteristic $p>0$ and $R:=k[x_0, \ldots ,
x_{n-1}]/(x_0^p, \ldots , x_{n-1}^p)$ (by Theorem
\ref{md11Feb06}.(2), $R$ is a differentiably simple Noetherian
commutative ring). Then
 $$\d :=\sum_{i=0}^{n-1} x^{[p^i-1]} \frac{\der }{\der x_i}\in \nsder
 (R)$$ where $x^{[p^i-1]}:=
 \prod_{\nu =0}^{i-1}\frac{x_\nu^{p-1}}{(p-1)!}$, $1\leq i \leq
 n-1$, and $x^{[0]}:=1$.
\end{lemma}

{\it Proof}. For each natural number  $j=0,1, \ldots , p^n-1$,
written $p$-adically, $j= \sum_\nu j_\nu p^\nu$, let
$x^{[j]}:=\prod_\nu \frac{x_\nu^{j_\nu}}{j_\nu !}$. By Proposition
\ref{ittp}, $\{ x^{[j]}, 0\leq j <p^n\}$ is the iterative sequence
with $x^{[0]}:=1$. By Corollary \ref{1ittp}, this iterative
sequence is a $\d$-descent. Since
$R=\bigoplus_{j=0}^{p^n-1}kx^{[j]}$ and $\d (x^{[j]})=x^{[j-1]}$
for all $j=0,1, \ldots , p^n-1$, the derivation $\d$ is simple and
nilpotent with $\d^{p^n}=0$. For each $i=0,1, \ldots , n-1$,
$\d^{p^i}(x^{[p^i]})=1$, and so $\d \in \nsder (R)$. $\Box $

%%%%%%%%%%%%%%%%%% SECTION 4 %%%%%%%%%%%%%%%%%%%%%%%%

\section{Simple derivations of differentiably simple Noetherian commutative rings}

In this section, we will see that for a differentially simple
Noetherian commutative ring $(R, \gm )$ there are  strong
connections between simple nilpotent derivations $\nsder (R)$,
coefficient fields for $R$, iterative descents, and the group
$\Aut (R)$ of ring automorphisms of $R$. Namely, there is a
canonical bijection $\nsder (R) \simeq \Aut (R)/ \Aut (R/ \gm )$
(Corollary \ref{a19Feb06}).

Recall that for a local commutative ring $(R, \gm )$,  a subfield
$k'$ of $R$ is called a coefficient field of $R$ if $R= k'+\gm $.
Let $\CF (R)$ be the set of all the coefficient fields of $R$. If
$k'$ is a coefficient  field of $R$, then $k'\simeq (k'+\gm )/ \gm
\simeq R/ \gm = k$, the residue field of $R$.

Let $k$ be  a field  of characteristic $p>0$. For  a
differentiably simple Noetherian commutative ring $(R, \gm )$ of
characteristic $p>0$ with $n:=\dim_k(\gm / \gm^2)\geq 1$ where
$k:= R/\gm$ (i.e.,  $R\simeq T_n:= k[x_0, \ldots , x_{n-1}]/
(x_0^p, \ldots , x_{n-1}^p)$), a subset $\{ x_0', \ldots ,
x_{n-1}' \}\subseteq \gm $ is called a {\em canonical set of
generators} for $R$ if there exists a coefficient field $k'$ of
$R$ such that $R= k'\langle x_0', \ldots , x_{n-1}' \rangle \simeq
k'[x_0', \ldots , x_{n-1}']/ (x_0'^p, \ldots , x_{n-1}'^p)$. Let
$\CC = \CC (R)$ be the set of all such $(k'; x_0', \ldots ,
x_{n-1}')$. Then $\CC =\bigcup_{k'\in \CF (R)}\CC (R, k')$ is
 a disjoint union of its subsets
 $$\CC (R, k')=\{ (k'; x_0', \ldots ,
 x_{n-1}')\, | \, (k'; x_0', \ldots , x_{n-1}')\in \CC (R)\}.$$

Let $\sder (R)$ be the set of all {\em simple} derivations of the
ring $R$. For a local commutative ring $(R, \gm )$,   let $\nsder
(R)$ be the set of all {\em nilpotent simple} derivations $\d$ of
the ring $R$ such that if $\d^{p^i}\neq 0$,  then
$\d^{p^i}(y_i)=1$ for some $y_i\in \gm$.

The next theorem gives  $(i)$ another description of all the
coefficient fields for differentiably simple Noetherian
commutative $\Fp$-algebra, $(ii)$ the canonical form of each
derivation $\d \in \nsder (R)$ (Theorem \ref{md11Feb06}.(2)) via
the unique iterative $\d$-descent. Recall that the set $\nsder
(R)$ is a nonempty set (Lemma \ref{NnsdR}).

\begin{theorem}\label{md11Feb06}%\marginpar{md11Feb06}
Let $(R, \gm )$ be a differentiably simple Noetherian commutative
$\Fp$-algebra with residue field $k=R/\gm$,  $n:=\dim_k(\gm /
\gm^2)\geq 1$, and $\d \in \nsder (R)$. Then
\begin{enumerate}
\item $\d^{p^{n-1}}\neq 0$ and $\d^{p^n}=0$. \item There exists  a
{\bf unique} iterative $\d$-descent $\{ x^{[i]}, 0\leq
 i <p^n\}$. Then   $x^{[i]}\in \gm$ for all  $i=1, \ldots , p^n-1$;  $R^\d \in \CF
 (R)$;
 $$R=\bigoplus_{i=0}^{p^n-1}R^\d x^{[i]}=R^\d \langle x_0, \ldots ,
 x_{n-1}\rangle \simeq R^\d [x_0, \ldots ,
 x_{n-1}]/(x_0^p, \ldots ,
 x_{n-1}^p), $$ where $x_j:= x^{[p^j]}$,  $x^{[i]}=
 \prod_{\nu =0}^t\frac{x_\nu^{i_\nu}}{i_\nu!}$,  $i=\sum_{\nu =0}^ti_\nu p^\nu$,
 $0\leq i_\nu<p$; and $$\d =\sum_{i=0}^{n-1} x^{[p^i-1]} \frac{\der }{\der x_i}\in \Der_{R^\d } (R).$$
 \item The map $\phi := \sum_{i=0}^{p^n-1} (-1)^i x^{[i]}\d^i
 :R\ra R$ is a {\bf projection} onto the direct summand $R^\d$ of $R$.
 In particular, $R^\d = \phi (R)$.
\end{enumerate}
\end{theorem}

{\it Proof}. Let $s$ be the natural number such that
$\d^{p^{s-1}}\neq 0$ but  $\d^{p^{s}}=0$. For each $j=0, \ldots ,
s-1$, let $\d_j:= \d^{p^j}$ and  fix an element $y_j\in \gm$ such
that $\d_j(y_j)=1$ (we can do this since $\d \in \nsder (R)$).
Clearly, $y_j^p=0$ for all $j$ (Lemma \ref{Rm30Jan06}). Let $\{
x^{[i]}, 0\leq i <p^s\}$ be an iterative $\d$-descent (Theorem
\ref{f18Jan06}), and let $x_j:=x^{[p^j]}$, for $j=0, \ldots ,
s-1$. Note that $R= N(\d , R)_{p^s-1}$ since $\d^{p^s}=0$. By
(\ref{NdApn}), %\marginpar{1NdApn}
\begin{equation}\label{1NdApn}
R= \bigoplus_{i=0}^{p^s-1} R^\d x^{[i]}\simeq R^\d [ x_0, \ldots ,
x_{s-1}]/(x_0^p, \ldots , x_{s-1}^p).
\end{equation}
Since the derivation $\d$ is simple, $R^\d$ is a field (if $\ga$
is an ideal of $R^\d$, then $\gb := \bigoplus_{i=0}^{p^s-1} \ga
x^{[i]}$ is an ideal of $R$ such that $\d (\gb ) \subseteq \gb $).
Then, by Theorem \ref{f18Jan06}.(3), the iterative $\d$-descent
$\{ x^{[i]}, 0\leq i <p^s\}$ is unique. By (\ref{1NdApn}), $\gm =
(x_0, \ldots , x_{s-1})$; hence $s=n$. Now, statements 1 and 2
follows from (\ref{1NdApn}).

The map $\phi$ is an endomorphism of the $R^\d$-module $R$.  For
each $j=1, \ldots , p^n-1$,
$$ \phi (x^{[j]})=\sum_{i=0}^j(-1)^ix^{[i]}x^{[j-i]}=
\sum_{i=0}^j(-1)^i{j\choose i}\cdot x^{[j]}= (1-1)^j\cdot
x^{[j]}=0\cdot x^{[j]}=0.$$ Hence $\phi$ is the projection onto
$R^\d$ in  view of the decomposition $R= \bigoplus_{i=0}^{p^n-1}
R^\d x^{[i]}$; see (\ref{1NdApn}). This proves statement $3$.
$\Box $

\begin{theorem}\label{b10Feb06}%\marginpar{b10Feb06}
Let $(R, \gm )$ be a differentiably simple Noetherian commutative
$\Fp$-algebra with residue field $k=R/\gm$ and $n:=\dim_k(\gm /
\gm^2)\geq 1$.
\begin{enumerate}
\item  Then the map
$$ g: \nsder (R)\ra \CC (R), \;\; \d \mapsto (\ker \, \d ;
x^{[p^0]}, x^{[p^1]}, \ldots , x^{[p^{n-1}]}),$$ is a bijection
(where $\{ x^{[i]}, 0\leq i <p^n\}$ is the iterative $\d$-descent
as in Theorem \ref{md11Feb06}) with the inverse map given by the
rule
$$ g^{-1}: \CC (R)\ra  \nsder (R), \;\; (k';x_0, \ldots , x_{n-1})\mapsto \sum_{i=0}^{n-1}
x^{[p^i-1]} \frac{\der }{\der x_i}\in \Der_{k'}(R),$$ where
$x^{[0]}:=1$ and $x^{[p^i-1]}:= \prod_{j=0}^{i-1}
\frac{x_j^{p-1}}{(p-1)!}$ (see Theorem \ref{md11Feb06}). \item For
each coefficient field $k'$ in $R$, the restriction $g_{k'}$ of
the map $g$ to the subset $\nsder_{k'}(R):= \{ \d \in \nsder (R)\,
| \, \d (k')=0\}$ of $\nsder(R)$  yields an isomorphism $g_{k'}:
\nsder_{k'}(R)\ra \CC (R, k')$.
\end{enumerate}
\end{theorem}

{\it Proof}. $1$. The map $g$ is well defined due to Theorem
\ref{md11Feb06}.(1,3), and, by Theorem \ref{md11Feb06}.(2), for
each derivation $\d\in \nsder (R)$, $\d
=\sum_{i=0}^{n-1}x^{[p^i-1]} \frac{\der }{\der x_i}\in \Der_{R^\d
}(R)$, where $\{ x^{[j]}, 0\leq j <p^n\}$ is the iterative
$\d$-descent, $x_i:= x^{[p^i]}$ and $x^{[p^i-1]}=
\prod_{j=0}^{i-1} \frac{(x^{[p^j]})^{p-1}}{(p-1)!}=
\prod_{j=0}^{i-1} \frac{x_j^{p-1}}{(p-1)!}$.

Conversely, for each $(k'; x_0, \ldots , x_{n-1})\in \CC (R)$, let
$x:= \{ x^{[i]}, 0\leq i <p^n\}$ be the corresponding iterative
sequence which exists, by Proposition \ref{ittp}, since
$x_0^p=\cdots = x_{n-1}^p=0$. The derivation $\d : =
\sum_{i=0}^{n-1} x^{[p^i-1]}\frac{\der}{\der x_i}\in \Der_{k'}(R)$
is nilpotent (by the very definition of $\d$), and $\d
(x^{[p^i]})= \d (x_i) = x^{[p^i-1]}$, $0\leq i \leq n-1$ (by the
very definition of $\d$). By Corollary \ref{2ittp} and Theorem
\ref{f18Jan06}.(3), $x$ is the iterative $\d$-descent. Hence $R$
is a $\d$-simple ring with $\ker (\d ) = k'$ (Lemma
\ref{Ky13Feb06}.(3)).

$2$. If $k'$ is a coefficient field of $R=k[x_0, \ldots ,
x_{n-1}]/(x_0^p, \ldots , x_{n-1}^p)$, then the fields $k'$ and
$k$ are isomorphic. Let us fix an isomorphism, say $\s : k\ra k'$.
Then $\s$ can be extended to an automorphism of the ring $R$ by
setting $\s (x_i)=x_i$ for all $i$ (since $R=k'+\gm = k'\langle
x_0, \ldots ,x_{n-1}\rangle \simeq k'[x_0, \ldots ,
x_{n-1}]/(x_0^p, \ldots , x_{n-1}^p)$). This implies that the set
$\nsder_{k'}(R)\neq \emptyset$ since $\nsder_k(R)\neq \emptyset$.
Now, statement $2$ follows from statement $1$.  $\Box $

For  the ring $R$ and its coefficient field $k'\in \CF (R)$, let
$\Aut (R)$ (resp. $\Aut_{k'}(R)$) be the group of all ring (resp.
$k'$-algebra) automorphisms of $R$. For the residue field $k:= R/
\gm$, $\Aut (k)$ is the group of its automorphisms.

Recall that an action of a group $G$ on a set $X$ is said to be
{\em fully faithful} if for some/each $x\in X$ the map $G\ra X$,
$g\mapsto gx$, is a bijection.

\begin{corollary}\label{a19Feb06}%\marginpar{a19Feb06}
Let $(R, \gm )$ be as in Theorem \ref{b10Feb06}. Then
\begin{enumerate}
\item  for each coefficient field $k'\in \CF (R)$, the action
 ${\rm Aut}_{k'}(R)\times \nsder_{k'} (R)\ra \nsder_{k'} (R)$
 which is given by the rule $(\s, \d )\mapsto \s \d \s^{-1}$
 is fully faithful.
\item  The action ${\rm Aut}(R)\times \nsder (R)\ra \nsder (R)$,
$(\s, \d )\mapsto \s \d \s^{-1}$, has a single orbit and,  for
each $\d \in \nsder (R)$, ${\rm Fix} (\d )\simeq {\rm Aut}(k)$,
and so $\nsder (R)\simeq {\rm Aut}(R)/{\rm Aut} (k)$ where ${\rm
Fix} (\d ):= \{ \s \in \Aut (R)\, | \, \s \d \s^{-1}= \d \}$.
\end{enumerate}
\end{corollary}

{\it Proof}. $1$. By the Theorem of Harper (Theorem
\ref{Harpbig}.(2)),  one can assume that $R=k'[x_0, \ldots ,
x_{n-1}]/ (x_0^p, \ldots , x_{n-1}^p)$. The natural action of the
group ${\rm Aut}_{k'}(R)$ on the set $\CC (R, k')$ is fully
faithful where $\s \cdot (k';x_0, \ldots , x_{n-1}):=(k'; \s
(x_0), \ldots , \s (x_{n-1}))$. The bijection $g_{k'}: \nsder_{k'}
(R) \ra \CC (R, k')$ {\em commutes} with the action of the group
${\rm Aut}_{k'}(R)$. Therefore, the action of ${\rm Aut}_{k'}(R)$
on $\nsder_{k'} (R)$ is fully faithful.

$2$. The action of $\Aut (R)$ on $\CC (R)$,  $\s \cdot (k;x_0,
\ldots , x_{n-1})):=(\s (k); \s (x_0), \ldots , \s (x_{n-1}))$,
has a single orbit. The bijection $g$ from Theorem
\ref{b10Feb06}.(1) commutes with the action of the group $\Aut
(R)$. Therefore, $\nsder (R)$ is an orbit of $\Aut (R)$.

One can assume that $R= k[x_0, \ldots , x_{n-1}]/ (x_0^p, \ldots ,
x_{n-1}^p)$ and $\d = \sum_{i=0}^{n-1} x^{[p^i-1]}
\frac{\der}{\der x_i}\in \Der_k(R)$ where $k = \ker \, \d$ (see
Theorem \ref{md11Feb06}). Consider the group monomorphism $\Aut
(k)\ra \Aut (R)$, $\tau \mapsto \tau$, given by the rule $\tau
(x_i)= x_i$ for all $i$. It remains to show that ${\rm Fix} (\d )
= \Aut (k)$. The inclusion $\Aut (k)\subseteq {\rm Fix } (\d )$ is
obvious. To prove the reverse inclusion we must show that if an
automorphism  $\s $ commutes with $\d$, then $\s (x_i) = x_i$ for
all $i$. Note that each $x_i':= \s (x_i)\in \gm $ for all $i$.
Since  $\d (x_0')= \d \s (x_0)=\s \d (x_0)= \s (1)=1$, we have
$x_0'= x_0+\l $ for some scalar  $\l \in k$. The scalar $\l $ must
be zero since $x_0'\in \gm$, and so $x_0'=x_0$. Suppose that
$x_0'=x_0, \ldots ,x_{i-1}'=x_{i-1}$ for some $i\geq 1$. Then
$x^{[p^i-1]}=\prod_{k=0}^{i-1} \frac{x_k^{p-1}}{(p-1)!}$ and so
$\s (x^{[p^i-1]}) = x^{[p^i-1]}$. Now, $\d (x_i')= \s \d
(x^{[p^i]})= \s (x^{[p^i-1]})= x^{[p^i-1]} = \d (x_i)$. This
yields  the equality  $x_i'= x_i+\l $ for some scalar $\l \in k$.
The scalar $\l $   is forced to be zero since $x_i', x_i\in \gm$.
By induction, we have $x_j'=x_j$ for all $j$. $\Box $

\begin{proposition}\label{b17Feb06}%\marginpar{b17Feb06}
Let $(R, \gm )$ be a differentiably simple Noetherian commutative
$\Fp$-algebra, $k:= R/ \gm$, $n:= \dim_k(\gm / \gm^2)\geq 1$.
Then, for each $\d \in \nsder (R)$, $R^\d \in \CF (R)$ and
$\Der_{R^\d } (R) =\bigoplus_{i=0}^{n-1} R\d^{p^i}$.
\end{proposition}

{\it Proof}. Let $\{ x^{[i]}, 0\leq i <p^n\}$ be the iterative
$\d$-descent and $k':= R^\d$ (see Theorem \ref{md11Feb06}.(2)). By
Theorem \ref{md11Feb06}, $R\simeq k'[x_0, \ldots ,
x_{n-1}]/(x_0^p, \ldots , x_{n-1}^p)$ where $x_i:= x^{[p^i]}$,
$\gm = (x_0, \ldots , x_{n-1})$, and $k'\in \CF (R)$. It is
obvious that $\Der_{k'}(R)=\bigoplus_{i=0}^{n-1} R\der_i$ where
$\der_i := \frac{\der}{\der x_i}\in \Der_{k'}(R)$ (since for each
$\d \in \Der_{k'}(R)$, $\d = \sum_{i=0}^{n-1} \d (x_i) \der _i$),
and that
$$\d^{p^i}(x_j)\equiv \d_{ij}\equiv \der_i(x_j)\mod \gm \;\; {\rm  for\;  all}\;\;  i,j=0,
\ldots , n-1.$$ Hence $\Der_{k'}(R) =\sum_{i=0}^{n-1} R\d^{p^i}
+\gm \Der_{k'}(R)$. By the Nakayama Lemma, $\Der_{k'}
(R)=\sum_{i=0}^{n-1} R\d^{p^i}$. It is obvious that this sum is a
direct one (if $\der := \l_s\d^{p^s}+\cdots + \l_t\d^{p^t}=0$ is a
nontrivial relation, i.e.,  $\l_s\neq 0$, then $0= \der (x_s) =
\l_s$, a contradiction).  $\Box $

$${\bf Acknowledgements}$$

The author would like to thank the referee of this paper for a
careful reading of the manuscript and valuable comments.

Department of Pure Mathematics

University of Sheffield

Hicks Building

Sheffield S3 7RH

UK

email: v.bavula@sheffield.ac.uk

\end{document}